\documentclass[a4paper,12pt]{article}
\usepackage{amsfonts,amsmath,amssymb}
\makeatletter
\newbox\bk@bxb
\newbox\bk@bxa
\newif\if@bkcont
\newcount\bk@lcnt

\def\breakboxskip{2pt}
\def\breakboxparindent{1.8em}

\def\breakbox{\vskip\breakboxskip\relax
\setbox\bk@bxb\vbox\bgroup
\advance\linewidth -2\fboxrule
\hsize\linewidth\@parboxrestore
\parindent\breakboxparindent\relax}

\def\bk@split{%
\@tempdimb\ht\bk@bxb 
\advance\@tempdimb\dp\bk@bxb
\setbox\bk@bxa\vsplit\bk@bxb to\z@ 
\setbox\bk@bxa\vbox{\unvbox\bk@bxa}
\setbox\@tempboxa\vbox{\copy\bk@bxa\copy\bk@bxb}
\advance\@tempdimb-\ht\@tempboxa
\advance\@tempdimb-\dp\@tempboxa}

\def\bk@addfsepht{%
\setbox\bk@bxa\vbox{\vskip\fboxsep\box\bk@bxa}}

\def\bk@addskipht{%
\setbox\bk@bxa\vbox{\vskip\@tempdimb\box\bk@bxa}}

\def\bk@addfsepdp{%
\@tempdima\dp\bk@bxa
\advance\@tempdima\fboxsep
\dp\bk@bxa\@tempdima}

\def\bk@addskipdp{%
\@tempdima\dp\bk@bxa
\advance\@tempdima\@tempdimb
\dp\bk@bxa\@tempdima}

\def\bk@line{%
\hbox to \linewidth{%
\hskip-2\fboxsep\vrule \@width\fboxrule\hskip.5\fboxsep\vrule \@width\fboxrule\hskip1.5\fboxsep
\box\bk@bxa\hfil
}}%

\def\endbreakbox{\egroup
\ifhmode\par\fi{\noindent\bk@lcnt\@ne
\@bkconttrue\baselineskip\z@\lineskiplimit\z@
\lineskip\z@\vfuzz\maxdimen
\bk@split\bk@addfsepht\bk@addskipdp
\ifvoid\bk@bxb 
\def\bk@fstln{\bk@addfsepdp
\hskip-\parindent\vbox{\llap{\raisebox{-2ex}{\rule{1.5\fboxsep}{\fboxrule}\hskip.5\fboxsep}}\bk@line\llap{\rule{1.5\fboxsep}{\fboxrule}\hskip.5\fboxsep}}}

\else 
\def\bk@fstln{\vbox{\llap{\raisebox{-2ex}{\rule{1.5\fboxsep}{\fboxrule}\hskip.5\fboxsep}}\bk@line}\hfil%
\advance\bk@lcnt\@ne
\loop
\bk@split\bk@addskipdp\leavevmode
\ifvoid\bk@bxb 
\@bkcontfalse\bk@addfsepdp
\vtop{\bk@line\noindent\hskip-2\fboxsep{\rule{1.5\fboxsep}{\fboxrule}}}%

\else 
\bk@line
\fi
\hfil\advance\bk@lcnt\@ne
\if@bkcont\repeat}%
\fi
\leavevmode\bk@fstln\par}\vskip\breakboxskip\relax}

\def\smp{\smallskip\par}
\def\un{{\bf 1}}
\def\zero{\{0\}}
\def\pf{\noindent{\bf Proof~:}\ }
\def\findemo{~\leaders\hbox to 1em{\hss\  \hss}\hfill~\raisebox{.5ex}{\framebox[1ex]{}}\smp}
\def\spn{\bigskip\par\noindent}
\def\mpn{\medskip\par\noindent}
\def\smpn{\smallskip\par\noindent}
\def\normal{\mathop{\trianglelefteq}}

\def\smp{\smallskip\par}
\def\smpn{\smallskip\par\noindent}

\def\mpoint{\;\;.}
\def\mvirg{\;\;,}

\def\Res{{\rm Res}}

\def\Inf{{\rm Inf}}
\def\Def{{\rm Def}}

\def\Hom{{\rm Hom}}

\def\Inf{{\rm Inf}}
\def\Im{{\rm Im}}

\def\Aut{{\rm Aut}}

\def\Id{{\rm Id}}

\def\op{^{op}}

\def\Z{\mathbb{Z}}
\def\N{\mathbb{N}}
\def\F{\mathbb{F}}
\def\Q{\mathbb{Q}}

\newcommand{\romain}[1]{\uppercase\expandafter{\romannumeral #1}}

\newcommand{\flh}[2]{\mathop{\hbox to 12mm{\rightarrowfill}}_{\displaystyle #2}^{\displaystyle #1}\limits}
\newcommand{\sflh}[2]{\mathop{\hbox to 12mm{\rightarrowfill}}_{\scriptstyle #2}^{\scriptstyle #1}\limits}

\newcommand{\gset}[1]{#1\hbox{-$\mathsf{set}$}}

\newcommand{\sur}[1]{\,\overline{\! #1}}
\newcommand{\sumb}[2]{\mathop{\sum}_{{\scriptstyle #1}\atop {\scriptstyle #2}}}

\newcommand{\sumc}[3]{\sum_{{\scriptstyle #1}\atop {{\scriptstyle #2}\atop {\scriptstyle #3}}}}

\def\op{^{op}}

\newcommand{\carre}[8]{\begin{array}{ccc}
#1&\mathop{\hbox to 12mm{\rightarrowfill}}^{\displaystyle{#2}}\limits&#3\\
\llap{$\displaystyle{#4}$}\left\downarrow\vbox to 6mm{}\right. & & \left\downarrow\vbox to 6mm{}\right.\rlap{$\displaystyle{#5}$}\\
#6&\mathop{\hbox to 12mm{\rightarrowfill}}_{\displaystyle #7}\limits&#8\\
\end{array}}

\newcommand{\carrem}[8]{\begin{array}{ccc}
#1&\mathop{\hbox to 12mm{\rightarrowfill}}^{\displaystyle #2}\limits&#3\\
\llap{$\displaystyle #4$}\left\uparrow\vbox to 6mm{}\right. & & \left\uparrow\vbox to 6mm{}\right.\rlap{$\displaystyle #5$}\\
#6&\mathop{\hbox to 12mm{\rightarrowfill}}_{\displaystyle #7}\limits&#8\\
\end{array}}

\newenvironment{enonce}[1]{\pagebreak[2]\refstepcounter{subsection}\refstepcounter{prop}\smpn{{\bf \thesection.\arabic{prop}.\ \ #1~:}}\begin{it} }{\end{it}\smp}
\newenvironment{enonce*}[1]{\pagebreak[2]\smpn{#1~:}\begin{it} }{\end{it}\smp}
\newcommand{\result}[1]{\begin{enonce}{#1}}
\def\fresult{\end{enonce}}
\newcommand{\npar}{\smallskip\par\noindent\pagebreak[2]\refstepcounter{subsection}\refstepcounter{prop}{\bf \thesection.\arabic{prop}.\ \ }}

\newenvironment{mth}[1]{\begin{breakbox}\begin{enonce}{#1}}{\end{enonce}\end{breakbox}}
\newenvironment{mth*}[1]{\begin{breakbox}\begin{enonce*}{#1}}{\end{enonce*}\end{breakbox}}
\newenvironment{rem}[1]{\refstepcounter{subsection}\refstepcounter{prop} \mpn{{\bf \thesection.\arabic{prop}.}\ \ \bf#1\ :}}{\smp}

\def\dom{\backslash}
\makeatletter
\renewenvironment{enumerate}{\ifnum \@enumdepth >3 \@toodeep\else
      \advance\@enumdepth \@ne
      \edef\@enumctr{enum\romannumeral\the\@enumdepth}\list
      {\csname label\@enumctr\endcsname}{\setlength{\topsep}{1ex}\setlength{\itemsep}{0pt}\usecounter
        {\@enumctr}\def\makelabel##1{\hss\llap{##1}}}\fi}{\endlist}
\renewenvironment{itemize}{\ifnum \@itemdepth >3 \@toodeep\else \advance\@itemdepth \@ne
\edef\@itemitem{labelitem\romannumeral\the\@itemdepth}%
\list{\csname\@itemitem\endcsname}{\setlength{\topsep}{1ex}\setlength{\itemsep}{0pt}\def\makelabel##1{\hss\llap{##1}}}\fi}
{\endlist}
\makeatother

\makeatletter
\def\@sect#1#2#3#4#5#6[#7]#8{\ifnum #2>\c@secnumdepth
    \let\@svsec\@empty\else
    \refstepcounter{#1}\edef\@svsec{\csname the#1\endcsname .\hskip .5em}\fi
    \@tempskipa #5\relax
     \ifdim \@tempskipa>\z@
       \begingroup #6\relax
         \@hangfrom{\hskip #3\relax\@svsec}{\interlinepenalty \@M #8\par}%
       \endgroup
      \csname #1mark\endcsname{#7}\addcontentsline
        {toc}{#1}{\ifnum #2>\c@secnumdepth \else
                     \protect\numberline{\csname the#1\endcsname}\fi
                   #7}\else
       \def\@svsechd{#6\hskip #3\relax  
                  \@svsec #8\csname #1mark\endcsname
                     {#7}\addcontentsline
                          {toc}{#1}{\ifnum #2>\c@secnumdepth \else
                            \protect\numberline{\csname the#1\endcsname}\fi
                      #7}}\fi
    \@xsect{#5}}
\def\section{\@startsection {section}{1}{\z@}{-3.5ex plus-1ex minus
    -.2ex}{2.3ex plus.2ex}{\reset@font\Large\bf}}  

\makeatother

\renewenvironment{equation}{\refstepcounter{subsection}\refstepcounter{prop}$$}{\leqno{\bf (\theprop)}$$}

\def\mar[#1]{\ar@{-}[#1]|-{\object@{<}}}
\def\marb[#1]{\ar@{-}[#1]|{\object+{  }}}

\usepackage[all]{xy}
\newcommand{\mor}[1]{#1\hbox{-}\mathsf{Mor}}
\newcommand{\galmor}[1]{#1\hbox{-}\mathsf{Mor}^{\rm Gal}}
\def\Gal{^{\rm Gal}}
\newcommand{\gsec}[3]{#3/#2\to#3/#1}
\newcommand{\gsect}[1]{\langle #1\rangle}
\newcommand{\gen}[1]{{<}#1{>}}
\def\Gringo{\Xi}

\newcommand{\comp}[1]{^{\;\widehat{}}_{#1}}
\newcommand{\compn}[1]{\raisebox{-.5ex}{$^{\;\widehat{\widehat{}}}$}_{\!\!#1}}
\newcommand{\compni}[1]{\raisebox{-.6ex}{\footnotesize$^{\;\widehat{\widehat{}}}$}_{\!\!#1}}
\newcommand{\plusn}[1]{\raisebox{.2ex}{$^{\widehat{\scriptscriptstyle +}}$}_{\!\!\!\!#1}}
\newcommand{\plus}[1]{^{+}_{#1}}
\newcommand{\lln}{{\preceq\hspace{-.2ex}\rule{.08ex}{1ex}}}
\newcommand{\Lln}{{\preceq\hspace{-.9ex}\rule{.08ex}{1.5ex}}}
\begin{document}
\centerline{\Large\bf The slice Burnside ring}\vspace{1ex}
\centerline{\Large\bf and the section Burnside ring}\vspace{1ex}
\centerline{\Large\bf of a finite group}\vspace{.5cm}\par
\centerline{Serge Bouc}\vspace{1.2cm}\par
{\footnotesize{\bf Abstract~:} This paper introduces two new Burnside rings for a finite group $G$, called {\em the slice Burnside ring} and {\em the section Burnside ring}. They are built as Grothendieck rings of the category of morphisms of $G$-sets, and of {\em Galois morphisms} of $G$-sets, respectively. The well known results on the usual Burnside ring, concerning ghost maps, primitive idempotents, and description of the prime spectrum, are extended to these rings. It is also shown that these two rings have a natural structure of Green biset functor. The functorial structure of unit groups of these rings is also discussed.
\vspace{1cm}\par
{\bf AMS Subject Classification~:} 19A22, 18F30.\par
{\bf Keywords~:} Burnside ring, Green biset functor, slice, section, unit group.
}
\section{Introduction}
This paper introduces two variations on Burnside rings of finite groups, called {\em the slice Burnside ring} and {\em the section Burnside ring}. Both of them are built as Grothendieck rings of some category of {\em morphisms} of finite $G$-sets, instead of the category of finite $G$-sets used to build the usual Burnside ring. The difference between these two new Burnside rings is that the slice Burnside ring is built from arbitrary morphisms of finite $G$-sets, whereas the section Burnside ring uses only {\em Galois morphisms} of finite $G$-sets.\par
It turns out that most of the well known properties of the Burnside ring (see e.g.\cite{handbook}) extend to the slice Burnside ring and to the section Burnside ring~: both are commutative rings, which are free of finite rank as $\Z$-module. There is an analogue of Burnside's theorem~: both of these rings embed in a product of copies of the integers, via a {\em ghost map}, and this map has a finite cokernel. After tensoring with $\Q$, both rings become split semisimple $\Q$-algebras, and explicit formulae for their primitive idempotents can be stated. The prime spectrum of both rings can also be described, and Dress's characterization of solvable groups in terms of the connectedness of the spectrum of the Burnside ring can be generalized as well. Finally, both constructions have a natural biset functor structure, for which they become {\em Green biset functors}.\par
A major exception in this list of generalizable properties concerns unit groups~: it can be shown that, unlike the case of the usual Burnside ring, the correspondence sending a finite group to the unit group of either the slice or the section Burnside ring cannot be endowed with a structure of biset functor. This is due to the lack of a suitable {\em tensor induction} for these rings.\par
The paper is divided in two parts, and an appendix~: the first part is devoted to the slice Burnside functor. It consists of Sections~\ref{morphisms} to~\ref{Prime spectrum}. Section~\ref{morphisms} recalls the basic definitions and properties on the category of morphisms of $G$-sets. Section~\ref{the slice Burnside functor} introduces the slice Burnside functor and its Green biset functor structure. Section~\ref{ghost} is devoted to the definition and main result on the ghost map. In Section~\ref{Slices and idempotents}, the explicit formulas for the primitive idempotents of the slice Burnside algebra over $\Q$ are stated. Section~\ref{ghost map image} gives a characterization of the image of the ghost map. Next Section~\ref{Prime spectrum} considers the prime spectrum of the slice Burnside ring. Finally, in Section~\ref{units}, it is shown that tom Dieck's theorem, building on Dress's characterization of solvable groups, can be extended to the slice Burnside ring~: namely, Feit-Thompson's theorem is equivalent to the fact that the only units in the slice Burnside ring of a group of odd order are $\pm 1$. But unlike the case of the usual Burnside ring, the unit group of the slice Burnside ring cannot be endowed with a structure of biset functor.\par
The second part of the paper is devoted to the section Burnside ring, and it is organized similarly~: Section~\ref{Galois morphisms} introduces Galois morphisms, and states their general properties. In particular, a left adjoint functor to the forgetful functor from Galois morphisms to arbitrary morphisms of $G$-sets is described, which attaches to any morphism of $G$-sets a canonical Galois morphism. Section~\ref{The section Burnside functor} considers the section Burnside functor, with its Green biset functor structure. Section~\ref{Galois ghost} deals with the ghost map for the section Burnside ring. Sections~\ref{Sections and idempotents} states the formulas for the primitive idempotents of the section Burnside $\Q$-algebra, and Section~\ref{Galois ghost map image} gives a characterization of the ghost map. Section~\ref{Galois prime spectrum} considers the prime spectrum of the section Burnside ring. In Section~\ref{Galois units}, the results about unit groups of Section~\ref{units} are extended to the unit group of the section Burnside ring. In particular, it is shown that this unit group cannot be endowed with a structure of biset functor.\par
The appendix deals with the functorial structure of the unit group of the slice Burnside ring and the section Burnside ring~: it is possible to define biset functor operations for these unit groups, but only for {\em left inert} bisets. This gives two interesting examples of somehow natural biset functors {\em without induction}. The last result of the appendix is the explicit computation of this unit group in the case of an abelian group.
\vspace{8ex}\par
\centerline{\Large\bf \romain{1} - The slice Burnside ring}
\section{Morphisms of $G$-sets}\label{morphisms}
\begin{mth}{Definition} Let $G$ be a group. If $f:X\to Y$ and $f':X'\to Y'$ are morphisms of $G$-sets, a {\em morphism} from $f$ to $f'$ is a pair of morphisms of $G$-sets $\alpha:X\to X'$ and $\beta: Y\to Y'$ such that the diagram
$$\xymatrix{
X\ar[r]^-{f}\ar[d]_{\alpha}&Y\ar[d]^-{\beta}\\
X'\ar[r]^-{f'}&Y'\\
}$$
is commutative.\par
Morphisms of morphisms of $G$-sets can be composed in the obvious way. This composition endows the class of morphisms of $G$-sets with a structure of category, denoted by $\mor{G}$.\par
\end{mth}
\begin{mth}{Proposition} \label{categorical product} The disjoint union of $G$-sets induces a coproduct
$$(X\stackrel{f}{\to}Y,X'\stackrel{f'}{\to}Y')\mapsto (X\sqcup X'\stackrel{f\sqcup f'}{\longrightarrow}Y\sqcup Y')$$
in the category $\mor{G}$.\par
Similarly, the direct product of $G$-sets, with diagonal $G$-action, induces a product
$$(X\stackrel{f}{\to}Y,X'\stackrel{f'}{\to}Y')\mapsto (X\times X'\stackrel{f\times f'}{\longrightarrow}Y\times Y')$$
in the category $\mor{G}$.\par
\end{mth}
\pf For any morphism of $G$-sets, the bijections
$$\Hom_{G\hbox{-}\mathsf{Set}}(X\sqcup X',A)\cong\Hom_{G\hbox{-}\mathsf{Set}}(X,A)\times \Hom_{G\hbox{-}\mathsf{Set}}(X',A)$$
induce obvious bijections between
$$\Hom_{\mor{G}}\big((X\sqcup X')\stackrel{f\sqcup f'}{\longrightarrow}(Y\sqcup Y'),A\stackrel{\alpha}{\to}B\big)$$
and
$$\Hom_{\mor{G}}(X\stackrel{f}{\to}Y,A\stackrel{\alpha}{\to}B)\times \Hom_{\mor{G}}(X'\stackrel{f'}{\to}Y',A\stackrel{\alpha}{\to}B)\mpoint$$
These bijections are obviously functorial in $\mor{G}$.\par
Similarly, the bijections
$$\Hom_{G\hbox{-}\mathsf{Set}}(A,X\times X')\cong\Hom_{G\hbox{-}\mathsf{Set}}(A,X)\times \Hom_{G\hbox{-}\mathsf{Set}}(A,X')$$
induce obvious bijections between
$$\Hom_{\mor{G}}\big(A\stackrel{\alpha}{\to}B,(X\times X')\stackrel{f\times f'}{\longrightarrow}(Y\times Y')\big)$$
and
$$\Hom_{\mor{G}}(A\stackrel{\alpha}{\to}B,X\stackrel{f}{\to}Y)\times \Hom_{\mor{G}}(A\stackrel{\alpha}{\to}B,X'\stackrel{f'}{\to}Y')\mpoint$$
These bijections are obviously functorial in $\mor{G}$.\findemo
\section{The slice Burnside functor} \label{the slice Burnside functor}
\begin{mth}{Definition and Notation} Let $G$ be a group. A {\em slice} of $G$ is a pair $(T,S)$ of subgroups of $G$ with $T\geq S$. A {\em section} of $G$ is a slice $(T,S)$ with $S\normal T$.\par
Let $\Pi(G)$ denote the set of slices of $G$, and $\Sigma(G)$ denote the set of sections of $G$.\par
When $(T,S)\in\Pi(G)$, denote by $G/S\to G/T$ the projection morphism.
\end{mth}
\vspace{1ex}
\begin{mth}{Definition} Let $G$ be a finite group. {\em The slice Burnside group $\Gringo(G)$ of $G$} is the quotient of the free abelian group on the set of isomorphism classes $[X\stackrel{f}{\longrightarrow}Y]$ of morphisms of finite $G$-sets, by the subgroup generated by elements of the form
$$[(X_1\sqcup X_2)\stackrel{f_1\sqcup f_2}{\longrightarrow}Y]-[X_1\stackrel{f_1}{\to}f(X_1)]-[X_2\stackrel{f_2}{\to}f(X_2)]\mvirg$$
whenever $X\stackrel{f}{\to}Y$ is a morphism of finite $G$-sets with a decomposition $X=X_1\sqcup X_2$ as a disjoint union of $G$-sets, where $f_1=f_{\mid X_1}$ and $f_2=f_{\mid X_2}$ .\par
When $f:X\to Y$ is a morphism of finite $G$-sets, let $\pi(f)$ denote the image in $\Gringo(G)$ of the isomorphism class of $f$. If $S\leq T$ are subgroups of $G$, set $\langle T,S\rangle_G=\pi(G/S\to G/T)$.
\end{mth}
\vspace{2ex}\pagebreak[3]
\begin{mth}{Lemma} \label{image}\begin{enumerate}
\item $\pi(\emptyset\to\emptyset)=0$.
\item Let $X\stackrel{f}{\to}Y$ be a morphism of finite $G$-sets. Then
$$\pi(X\stackrel{f}{\to}Y)=\pi\big(X\stackrel{f}{\to}f(X)\big)\mpoint$$
\item Let $X\stackrel{f}{\to}Y$ and $X'\stackrel{f'}{\to}Y'$ be morphisms of finite $G$-sets. Then
$$\pi\big((X\sqcup X')\stackrel{f\sqcup f'}{\longrightarrow}(Y\sqcup Y')\big)=\pi(X\stackrel{f}{\to}Y)+\pi(X'\stackrel{f'}{\to}Y')\mpoint$$
\end{enumerate}
\end{mth}
\pf For Assertion 1, set $e= \pi(\emptyset\to\emptyset)$. Since the morphism $\emptyset \sqcup\emptyset\to\emptyset$ is isomorphic to $\emptyset\to\emptyset$, it follows that $e+e=e$, hence $e=0$.\par
For Assertion 2, writing $X=X\sqcup \emptyset$ gives 
$$\pi\big(X\stackrel{f}{\to} Y\big)=\pi\big(X\stackrel{f}{\to} f(X)\big)+\pi(\emptyset\to\emptyset)=\pi\big(X\stackrel{f}{\to} f(X)\big)\mpoint$$
For Assertion 3, 
\begin{eqnarray*}
\pi\big((X\sqcup X')\stackrel{f\sqcup f'}{\longrightarrow}(Y\sqcup Y')\big)&=&\pi\big(X\stackrel{f}{\to}(Y\sqcup Y')\big)+\pi\big(X'\stackrel{f'}{\to}(Y\sqcup Y')\big)\\
&=&\pi\big(X\stackrel{f}{\to}f(X)\big)+\pi\big(X'\stackrel{f'}{\to}f'(X')\big)\\
&=&\pi\big(X\stackrel{f}{\to}Y\big)+\pi\big(X'\stackrel{f'}{\to}Y'\big)\mvirg
\end{eqnarray*}
where the first equality follows from the defining relations of $\Gringo(G)$, and the other ones from Assertion~2.\findemo
\begin{mth}{Lemma} \label{decomposition}Let $f:X\to Y$ be a morphism of finite $G$-sets. Then in the group $\Gringo(G)$
\begin{eqnarray*}
\pi(X\stackrel{f}{\to}Y)&=&\sum_{x\in [G\dom X]}\langle G_{f(x)},G_x\rangle_G\\
\end{eqnarray*}
\end{mth}
\pf Indeed $X\cong \mathop{\sqcup}_{x\in[G\dom X]}\limits G/G_x$, and the image $f(G\cdot x)$ of the $G$-orbit of $x$ is equal to the $G$-orbit of $f(x)$. Moreover the morphisms $f_{\mid G\cdot x}:G\cdot x\to G\cdot f(x)$ and $\gsec{G_{f(x)}}{G_x}{G}$ are isomorphic. The claimed formula follows. 
\findemo
\begin{mth}{Corollary} \label{generators}The group $\Gringo(G)$ is generated by the elements $\langle T,S\rangle_G$, where $(T,S)$ runs through a set $[\Pi(G)]$ of representatives of conjugacy classes of slices of $G$.
\end{mth}
\pf Indeed, the morphisms $G/{^gS}\to G/{^gT}$ and $G/S\to G/T$ are isomorphic, for any $g\in G$, and any slice $(T,S)$ of $G$.\findemo
\begin{rem}{Remark} It will be shown in Theorem~\ref{Gringo basis} that this generating set is actually a basis of $\Gringo(G)$.
\end{rem}
\begin{mth}{Proposition} The product of morphisms induces a commutative unital ring structure on $\Gringo(G)$. The identity element for multiplication is the image of the class $[\bullet\to\bullet]$, where $\bullet$ denotes a $G$-set of cardinality 1.
\end{mth}
\pf If we can show that the product of morphisms induces a well defined bilinear product $\Gringo(G)\times\Gringo(G)\to \Gringo(G)$, it will be clear that this product is associative, commutative, and admits $[\bullet\to\bullet]$ as an identity element. Hence the only point to check is that the product preserves the defining relations of $\Gringo(G)$. This is clear, since if $g:Z\to T$ is a morphism of finite $G$-sets, and if $X_1\sqcup X_2\stackrel{f_1\sqcup f_2}{\longrightarrow}Y$ is a morphism, setting $X=X_1\sqcup X_2$, the domain of the morphism
$$
\xymatrix{
h:Z\times X\ar[rr]^-{g\times (f_1\sqcup f_2)}&&T\times Y
}
$$
has a disjoint union decomposition $Z\times X=(Z\times X_1)\sqcup (Z\times X_2)$, and moreover the restriction of $g\times f$ to $Z\times X_1$ is $g\times f_1$. Thus
\begin{eqnarray*}
\pi(h)\!\!&\!\!=\!\!&\!\!\pi\Big((Z\times X_1)\stackrel{g\times f_1}{\longrightarrow}\big(g(Z)\times f_1(X_1)\big)\Big)+\pi\Big((Z\times X_2)\stackrel{g\times f_2}{\longrightarrow}\big(g(Z)\times f_2(X_2)\big)\Big)\\
\!\!&\!\!=\!\!&\!\!\pi\big((Z\times X_1)\stackrel{g\times f_1}{\longrightarrow}(T\times f_1(X_1))\big)+\pi\big((Z\times X_2)\stackrel{g\times f_2}{\longrightarrow}(T\times f_2(X_2))\big)
\end{eqnarray*}
where the last equality follows from Lemma~\ref{image}.\findemo
\begin{mth}{Proposition} \label{product formula}Let $(T,S)$ and $(Y,X)$ be slices of $G$. Then in $\Gringo(G)$
$$\gsect{T,S}_G\gsect{Y,X}_G=\sum_{g\in[S\dom G/X]}\gsect{T\cap {^gY},S\cap{^gX}}_G\mpoint$$
\end{mth}
\pf Indeed
$$(G/S)\times (G/X)\cong \mathop{\sqcup}_{g\in [S\dom G/X]}G/(S\cap {^gX})\mvirg$$
via the map (from right to left) sending $u(S\cap{^gX})$ to $(uS,ugX)$, for $u\in G$. The image of $(S,gX)$ by the map $(G/S)\times (G/X)\to (G/T)\times (G/Y)$ is the pair $(T,gY)$, whose stabilizer in $G$ is $T\cap{^gY}$. The result now follows from Lemma~\ref{decomposition}.\findemo
\pagebreak[3]
\begin{mth}{Theorem} \label{Green biset functor}\begin{enumerate}
\item Let $G$ and $H$ be finite groups, and let $U$ be a finite $(H,G)$-biset. The functor
$$(X\stackrel{f}{\to}Y)\mapsto (U\times_GX\stackrel{U\times_Gf}{\longrightarrow}U\times_GY)$$
from $\mor{G}$ to $\mor{H}$ induces a group homomorphism 
$$\Gringo(U):\Gringo(G)\to\Gringo(H)\mpoint$$
\item The correspondence $G\mapsto \Gringo(G)$ is a Green biset functor.
\end{enumerate}
\end{mth}
\pf For Assertion~1, the only thing to check is that the defining relations of $\Gringo(G)$ are mapped to relations in $\Gringo(H)$. But if 
$$X_1\sqcup X_2\stackrel{f_1\sqcup f_2}{\longrightarrow}Y$$
is a morphism of finite $G$-sets, then
$$U\times_G(X_1\sqcup X_2)\cong (U\times_GX_1)\sqcup(U\times_GX_2)\mpoint$$
Moreover the image of the map $U\times_Gf_1$ is equal to $U\times_Gf_1(X_1)$. It follows that the relation
$$[X_1\sqcup X_2\stackrel{f_1\sqcup f_2}{\longrightarrow}Y]-[X_1\stackrel{f_1}{\to}f_1(X_1)]-[X_2\stackrel{f_2}{\to}f_2(X_2)]$$
in $\Gringo(G)$ is mapped to the relation
$$[UX_1\sqcup UX_2\stackrel{Uf_1\sqcup Uf_2}{\longrightarrow}UY]-[UX_1\stackrel{Uf_1}{\to}(Uf_1)(UX_1)]-[UX_2\stackrel{Uf_2}{\to}(Uf_2)(UX_2)]\mvirg$$
where $U\times_G$ is abbreviated to $U$.\par
It is now clear that the correspondence sending a finite group $G$ to $\Gringo(G)$ and a finite $(H,G)$-biset $U$ to $\Gringo(U)$ endows $\Gringo$ with a structure of biset functor (see \cite{bisetfunctors}).\par
Moreover if $G$, $G'$ are finite groups, if $f:X\to Y$ is a morphism of finite $G$-sets and $f':X'\to Y'$ is a morphism of finite $G'$-sets,
then 
$$f\times f':X\times X'\to Y\times Y'$$
is a morphism of $G\times G'$-sets. This induces a product 
$$\Xi(G)\times \Xi(G')\to \Xi(G\times G')\mvirg$$
which is associative in the obvious sense. Moreover, the morphisms $\bullet\to\bullet$ of $\un$-sets is obviously an identity element for this product, up to identification $G\times\un=G$. \par
Finally, if $G$, $G'$, $H$, $H'$ are finite groups, if $U$ is a finite $(H,G)$-biset, if $U'$ is a finite $(H',G')$-biset, it is clear that the morphisms
$$(U\times U')\times_{G\times G'}(f\times f')$$
and 
$$(U\times_Gf)\times (U'\times_{G'}f')$$
are isomorphic morphisms of $(H\times H')$-sets.  Thus $\Gringo$ is a Green biset functor (see \cite{bisetfunctors} Section 8.5).\findemo
\begin{mth}{Proposition} Let $G$ and $H$ be finite groups, and $U$ be a finite $(H,G)$-biset. If $(T,S)\in\Pi(G)$, then
$$U\times_G\gsect{T,S}_G=\sum_{u\in [H\dom U/S]}\gsect{{^uT},{^uS}}_H\mpoint$$
(where $^uX=\{h\in H\mid \exists x\in X,\;hu=ux\}$, for $X\leq G$).
\end{mth}
\pf Indeed $U\times_G(G/S)\cong U/S$, and the stabilizer in $H$ of $uS$ is equal to ${^uS}$. The result follows from Lemma~\ref{decomposition}.\findemo
The following proposition clarifies the links existing between the usual Burnside functor and the slice Burnside functor~:
\begin{mth}{Proposition} \label{from B to Xi}\begin{enumerate}
\item Let $G$ be a finite group. The correspondence sending the morphism $X\stackrel{f}{\to}Y$ of finite $G$-sets to the $G$-set $X$ induces a unital ring homomorphism $s_G$ from $\Xi(G)$ to the Burnside ring $B(G)$.
\item The correspondence sending the finite $G$-set $X$ to the identity morphism of $X$ induces a unital ring homomorphism $i_G:B(G)\to\Xi(G)$, such that $s_G\circ i_G=\Id_{B(G)}$.
\item As $G$ varies, the morphisms $s_G$ and $i_G$ define morphisms of Green biset functors $s:\Xi\to B$ and $i:B\to \Xi$, such that $s\circ i=\Id_B$. In particular $i$ is injective, and $s$ is surjective.
\end{enumerate}
\end{mth}
\pf This is straightforward, from the definitions.\findemo
\section{Slices and ghost map}\label{ghost}
\vspace{-2ex}
\begin{mth}{Notation} Let $S\leq T$ be subgroups of $G$. If $X\stackrel{f}{\to}Y$ is a morphism of finite $G$-sets, set
$$\phi_{T,S}(X\stackrel{f}{\to}Y)=|\Hom_{\mor{G}}\big(\gsec{T}{S}{G},X\stackrel{f}{\to}Y)|\mpoint$$
\end{mth}
\begin{mth}{Notation} Define a relation $\preceq$ on the set $\Pi(G)$ by
$$(T,S)\preceq (Y,X)\Leftrightarrow (T\leq Y\;\;\hbox{and}\;\;S\leq{X})\mpoint$$
\vspace{-2ex}
\end{mth}
The relation $\preceq$ is an order relation on $\Pi(G)$.
\begin{mth}{Lemma} \label{phiTS}With this notation
$$\phi_{T,S}(X\stackrel{f}{\to}Y)=|f^{-1}(Y^T)^S|\mpoint$$
In particular, for any $A\leq B\leq G$,
$$\phi_{T,S}(G/A\to G/B)=|\{g\in G/A\mid (T^g,S^g)\preceq (B,A)\}|\mpoint$$
\end{mth}
\pf The morphisms of $G$-sets from $G/S$ to $X$ are in one to one correspondence with the set $X^S$ of fixed points of $S$ on $X$~: the morphism associated to $x\in X^S$ is defined by $gS\mapsto gx$. Similarly, the homomorphisms of $G$-sets from $G/T$ to $Y$ are in one to one correspondence with $Y^T$. Hence the set
$$\Hom_{\mor{G}}(G/S\to G/T,X\stackrel{f}{\to}Y)$$
is in one to one correspondence with the set of pairs $(x,y)\in X^S\times Y^T$ such that $f(x)=y$, i.e. with the elements $x$ of $f^{-1}(Y^T)^S$.\findemo
\begin{mth}{Corollary} \label{modp}Let $(T,S)\in\Pi(G)$ and $p$ be a prime number. If $P$ is a $p$-subgroup of $N_G(T,S)$, then
$$\phi_{T,S}\equiv\phi_{PT,PS}\;({\rm mod.}\;p)\mpoint$$
\end{mth}
\pf Indeed for any morphism of finite $G$-sets $X\stackrel{f}{\to}Y$, the set $f^{-1}(Y^T)^S$ is invariant by $N_G(T,S)$, thus, as $P$ is a $p$-group,  
$$|f^{-1}(Y^T)^S|\equiv |f^{-1}(Y^T)^{PS}|\;({\rm mod.}\;p)\mvirg$$
and moreover $f^{-1}(Y^T)^{PS}=f^{-1}(Y^{PT})^{PS}$.\findemo
\vspace{-1ex}
\begin{mth}{Proposition} Let $S\leq T$ be subgroups of $G$. Then the map $\phi_{T,S}$ induces a ring homomorphism $\Gringo(G)\to \Z$, still denoted by $\phi_{T,S}$.
\end{mth}
\pf Since the product of morphisms
$$(X\stackrel{f}{\to}Y,X'\stackrel{f'}{\to}Y') \mapsto \big((X\times X')\stackrel{f\times f'}{\longrightarrow}(Y\times Y')\big)$$
is a product in the category $\mor{G}$, it follows that 
$$\phi_{T,S}\big((X\times X')\stackrel{f\times f'}{\longrightarrow}(Y\times Y')\big)=\phi_{T,S}(X\stackrel{f}{\to}Y)\phi_{T,S}(X'\stackrel{f'}{\to}Y')\mpoint$$
Also $\phi_{T,S}(\bullet\to\bullet)=1$. The only thing to check is that $\phi_{T,S}$ induces a well defined map $\Gringo(G)\to \Z$, i.e. that the defining relations of $\Gringo(G)$ are mapped to 0 by $\phi_{T,S}$. First, by Lemma~\ref{phiTS}, for any morphism of $G$-sets $f:X\to Y$
\vspace{-1ex}
$$\phi_{T,S}(X\stackrel{f}{\to}Y)=\phi_{T,S}\big(X\stackrel{f}{\to}f(X)\big)\mpoint$$
\vspace{-2ex}
Now
\begin{eqnarray*}
\phi_{T,S}\big((X_1\sqcup X_2)\stackrel{f_1\sqcup f_2}{\longrightarrow}Y\big)\!&\!=\!&\!|(f_1\sqcup f_2)^{-1}(Y^T)^S|\\
\!&\!=\!&\!|X_1\cap (f_1\sqcup f_2)^{-1}(Y^T)^S|+|X_2\cap (f_1\sqcup f_2)^{-1}(Y^T)^S|\\
\!&\!=\!&\!|f_1^{-1}(Y^T)^S|+|f_2^{-1}(Y^T)^S|\\
\!&\!=\!&\!\phi_{T,S}(X_1\stackrel{f_1}{\to}Y)+\phi_{T,S}(X_2\stackrel{f_2}{\to}Y)\\
\!&\!=\!&\!\phi_{T,S}\big(X_1\stackrel{f_1}{\to}f(X_1)\big)+\phi_{T,S}\big(X_2\stackrel{f_2}{\to}f(X_2)\big)\mpoint
\end{eqnarray*}
This completes the proof.\findemo
\vspace{-2ex}
\begin{mth}{Theorem} \label{Gringo basis}The group $\Gringo(G)$ is a free abelian group, with basis the set of elements $\gsect{T,S}_G$, where $(T,S)$ runs through a set $[\Pi(G)]$ of representatives of conjugacy classes of slices of $G$. Moreover, the map (called the {\em ghost map} for $\Gringo(G)$)
$$\Phi=\prod_{(T,S)\in[\Pi(G)]}\phi_{T,S}:\Gringo(G)\to \prod_{(T,S)\in[\Pi(G)]}\Z$$
is an injective ring homomorphism, with finite cokernel as morphism of abelian groups.
\end{mth}
\pf By Lemma~\ref{generators}, the elements $\langle T,S\rangle_G$, for $(T,S)\in[\Pi(G)]$, generate $\Gringo(G)$. Suppose that there is a non zero linear combination in the kernel of $\Phi$
$$\Lambda=\sum_{(T,S)\in[\Pi(G)]}\lambda_{T,S}\,\langle T,S\rangle_G$$
with integer coefficients $\lambda_{T,S}\in\Z$, for $(T,S)\in[\Pi(G)]$. Extend $\lambda$ to a function $\Pi(G)\to \Z$, constant on conjugacy classes. Let $(Y,X)$ be an element of $\Pi(G)$, maximal for the relation $\preceq$, such that $\lambda_{Y,X}\neq 0$. Then since by Lemma~\ref{phiTS}
$$\phi_{Y,X}(G/S\to G/T)=|\{g\in G/S\mid (Y^g,X^g)\preceq (T,S)\}|\mvirg$$
it follows that
\begin{eqnarray*}
\phi_{Y,X}(\Lambda)&=&\sum_{(T,S)\in[\Pi(G)]}\lambda_{T,S}\phi_{Y,X}(G/S\to G/T)\\
&=&\lambda_{Y,X}\phi_{Y,X}(G/X\to G/Y)\\
&=&\lambda_{Y,X}|N_G(X,Y)/X|=0\mpoint
\end{eqnarray*}
Hence $\lambda_{Y,X}=0$, and this contradiction shows that $\Phi$ is injective. In particular, the elements $\langle T,S\rangle_G$, for $(T,S)\in[\Pi(G)]$, form a $\Z$-basis of $\Gringo(G)$. Thus $\Phi$ is an injective morphism between free abelian groups with the same finite rank, hence it has finite cokernel.\findemo
\begin{mth}{Corollary} \label{Gringo semisimple}Set $\Q\Gringo(G)=\Q\otimes_\Z\Gringo(G)$, and $\Q\Phi=\Q\otimes_\Z\Phi$. Then
$$\Q\Phi:\Q\Gringo(G)\to\prod_{(T,S)\in[\Pi(G)]}\Q$$
is an isomorphism of $\Q$-algebras.
\end{mth}
\section{Slices and idempotents}\label{Slices and idempotents}
By Corollary~\ref{Gringo semisimple}, the commutative $\Q$-algebra $\Q\Gringo(G)$ is split semisimple. Its primitive idempotents are indexed by slices of $G$, up to conjugation~: they are the inverse images under $\Q\Phi$ of the primitive idempotents of the algebra $\mathop{\prod}_{(T,S)\in[\Pi(G)]}\limits\Q$~:
\begin{mth}{Notation}
If $(T,S)\in\Pi(G)$, denote by $\xi_{T,S}^G$ the unique element of $\Q\otimes_\Z\Gringo(G)$ such that
$$\forall (Y,X)\in\Pi(G),\;\;\Q\phi_{Y,X}(\xi_{T,S}^G)=\left\{\begin{array}{cl}1&\hbox{if}\;(Y,X)=_G(T,S)\\0&\hbox{otherwise}\end{array}\right.$$
The set of elements $\xi_{T,S}^G$, for $(T,S)\in[\Pi(G)]$, is the set of primitive idempotents of $\Q\Gringo(G)$.
\end{mth}
\begin{mth}{Theorem} \label{slices and idempotents}Let $(T,S)\in \Pi(G)$. Then
$$\xi_{T,S}^G=\frac{1}{|N_G(T,S)|}\sum_{(V,U)\preceq (T,S)}|U|\mu_\Pi\big((V,U),(T,S)\big)\gsect{V,U}_G\mvirg$$
where $\mu_\Pi$ is the M\"obius function of the poset $(\Pi(G),\preceq)$.
\end{mth}
\pf
Denote by $\sur{\Pi}(G)$ the set of orbits of $G$ for its conjugacy action on $\Pi(G)$. Thus $\sur{\Pi}(G)$ is in one to one correspondence with $[\Pi(G)]$, and the map $\Q \Phi$ can also be viewed a $\Q$-algebra isomorphism from $\Q \Xi(G)$ to $\Q^{\sur{\Pi}(G)}$. The $\Q$-vector space $\Q\Xi(G)$ has a basis consisting of the elements $\gsect{V,U}_G$, for $(V,U)\in[\Pi(G)]$. Let $\Q^{\Pi(G)}$ denote the $\Q$-vector space with basis $\Pi(G)$, and let $p:\Q^{\Pi(G)}\to \Q\Xi(G)$ denote the $\Q$-linear map sending $(V,U)\in\Pi(G)$ to $\gsect{V,U}_G$. \par
Let $\beta_{T,S}^G$ denote the vector of the canonical basis of $\Q^{\sur{\Pi}(G)}$ indexed by the $G$-orbit of $(T,S)\in\Pi(G)$, and let $q:\Q^{\Pi(G)}\to \Q^{\sur{\Pi}(G)}$ denote the $\Q$-linear map sending $(T,S)\in\Pi(G)$ to $\beta_{T,S}^G$.\par
With this notation, Lemma~\ref{phiTS} shows that for any $(V,U)\in\Pi(G)$
\begin{eqnarray*}
\Phi\circ p\big((V,U)\big)&=&\sum_{(T,S)\in[\Pi(G)]}\phi_{T,S}(\gsect{V,U}_G)\beta_{T,S}^G\\
&=&\sum_{(T,S)\in[\Pi(G)]}\frac{1}{|U|}\big|\{g\in G\mid (T,S)\preceq ({^gV},{^gU})\}\big|\,\beta_{T,S}^G\\
&=&\sum_{(T,S)\in \Pi(G)}\frac{|N_G(T,S)|}{|G||U|}\big|\{g\in G\mid (T^g,S^g)\preceq ({V},{U})\}\big|\,\beta_{T,S}^G\\
&=&\sumc{(T,S)\in \Pi(G)}{g\in G}{(T^g,S^g)\preceq (V,U)}\frac{|N_G(T^g,S^g)|}{|G||U|}\,\beta_{T^g,S^g}^G\\
\end{eqnarray*}
Thus
$$\Phi\circ p\big((V,U)\big)=\sumb{(T,S)\in \Pi(G)}{(T,S)\preceq (V,U)}\frac{|N_G(T,S)|}{|U|}\,\beta_{T,S}^G$$
This shows that if $\widetilde{\Phi}$ is the $\Q$-endomorphism of $\Q^{\Pi(G)}$ defined by
$$\widetilde{\Phi}\big((V,U)\big)=\sumb{(T,S)\in \Pi(G)}{(T,S)\preceq (V,U)}\frac{|N_G(T,S)|}{|U|}(T,S)\mvirg$$
then $\Phi\circ p=q\circ \widetilde{\Phi}$. The matrix of the map $\widetilde{\Phi}$ is equal to the product $E\cdot J\cdot D$, where $D$ is a diagonal matrix with diagonal coefficients $(|N_G(T,S)|)_{(T,S)\in\Pi(G)}$, where $E$ is a diagonal matrix with diagonal coefficients $(\frac{1}{|U|})_{(V,U)\in\Pi(G)}$, and $J$ is the incidence matrix of the order relation $\preceq$ on $\Pi(G)$. It follows that $\widetilde{\Phi}$ is invertible, with inverse equal to $D^{-1}\cdot J^{-1}\cdot E^{-1}$. Now the entries of the matrix $J^{-1}$ are precisely the values of the M\"obius function $\mu_\Pi$ of the poset $\big(\Pi(G),\preceq\big)$. It follows that
\begin{eqnarray*}
\xi_{T,S}^G&=&\Phi^{-1}(\beta_{T,S}^G)\\
&=&\Phi^{-1}\circ q\big((T,S)\big)=p\circ\widetilde{\Phi}^{-1}\big((T,S)\big)\\
&=&\frac{1}{|N_G(T,S)|}\sumb{(V,U)\in\Pi(G)}{(V,U)\preceq (T,S)}|U|\mu_\Pi\big((V,U),(T,S)\big)\,\gsect{V,U}_G\mvirg
\end{eqnarray*}
which completes the proof.\findemo
\begin{mth}{Proposition} \label{poset of pairs}Let $(\mathcal{X},\leq)$ be a finite poset. Let $\Pi(\mathcal{X})$ denote the set of pairs $(y,x)$ of elements of $\mathcal{X}$ such that $x\leq y$. Define a partial order~$\preceq$ on $\Pi(\mathcal{X})$  by
$$\forall (y,x),\,(t,z)\in\Pi(\mathcal{X}),\;\;(y,x)\preceq (t,z)\Leftrightarrow y\leq t\;\hbox{and}\;x\leq z\mpoint$$
Then the M\"obius function $\mu_\Pi$ of the poset $(\Pi(\mathcal{X}),\preceq)$ can be computed as follows, for any $(y,x),\,(t,z)\in\Pi(\mathcal{X})$~:
\begin{equation}\label{mobius paires}\mu_\Pi\big((y,x),(t,z)\big)=\left\{\begin{array}{cl}\mu_\mathcal{X}(x,z)\mu_\mathcal{X}(y,t)&\hbox{if}\;\;x\leq z\leq y\leq t\\
0&\hbox{otherwise}\end{array}\right.\mvirg
\end{equation}
where $\mu_\mathcal{X}$ is the M\"obius function of the poset $(\mathcal{X},\leq)$.
\end{mth}
\pf Let $m\big((y,x),(t,z)\big)$ denote the expression defined by the right hand side of Equation~\ref{mobius paires}. Then if $(y,x)\preceq (t,z)$, i.e. if $y\leq t$ and $x\leq z$,
\begin{eqnarray*}
\sumb{(v,u)\in\Pi(\mathcal{X})}{\rule{0ex}{1.2ex}(y,x)\preceq (v,u)\preceq (t,z)}m\big((y,x),(v,u)\big)&=&\sumc{y\leq v\leq t}{\rule{0ex}{1.2ex}x\leq u\leq z}{\rule{0ex}{1.2ex}u\leq v,\,u\leq y}\mu_X(x,u)\mu_X(y,v)\\
&=&\sumc{y\leq v\leq t}{\rule{0ex}{1.2ex}x\leq u\leq z}{\rule{0ex}{1.2ex}u\leq y}\mu_X(x,u)\mu_X(y,v)\\
&=&\big(\sum_{y\leq v\leq t}\mu_\mathcal{X}(y,v)\big)\big(\sumb{x\leq u\leq z}{\rule{0ex}{1.2ex}u\leq y}\mu_\mathcal{X}(x,u)\big)\mpoint
\end{eqnarray*}
The first factor $\sum_{y\leq v\leq t}\limits\mu_\mathcal{X}(y,v)$ is equal to 0 if $y\neq t$, and to 1 if $y=t$. In this case, the second factor $\sumb{x\leq u\leq z}{\rule{0ex}{1.2ex}u\leq y}\limits\mu_\mathcal{X}(x,u)$ is equal to $\sum_{x\leq u\leq z}\limits\mu_\mathcal{X}(x,u)$. This is equal to zero if $x\neq z$, and to 1 if $x=z$. \par
It follows that $\sumb{(v,u)\in\Pi(\mathcal{X})}{\rule{0ex}{1.2ex}(y,x)\preceq (v,u)\preceq (t,z)}\limits m\big((y,x),(v,u)\big)$ is equal to $0$ if $(y,x)\neq (t,z)$ and to 1 otherwise. The proposition follows.\findemo
Applying Proposition~\ref{poset of pairs} to the poset of subgroups of $G$, ordered by inclusion of subgroups, gives the following~:
\begin{mth}{Corollary} \label{idempotents bis}Let $(V,U)$ and $(T,S)$ be slices of $G$. Then
\begin{equation}\label{poset Pi}\mu_\Pi\big((V,U),(T,S)\big)=\left\{\begin{array}{cl}\mu(U,S)\mu(V,T)&\hbox{if}\;U\leq S\leq V\leq T\\0&\hbox{otherwise}\end{array}\right.\mvirg
\end{equation}
where $\mu$ is the M\"obius function of the poset of subgroups of $G$.  In particular in $\Gringo(G)$
$$\xi_{T,S}^G=\frac{1}{|N_G(T,S)|}\sum_{U\leq S\leq V\leq T}|U|\mu(U,S)\mu(V,T)\,\gsect{V,U}_G\mpoint$$
\end{mth}
\section{The image of the ghost map}\label{ghost map image}
The following characterization of the image of the ghost map is the analogue for the slice Burnside ring of a theorem of Dress (\cite{dressresoluble}) on the ordinary Burnside ring ~:
\begin{mth}{Theorem} \label{ghost image} Let $G$ be a finite group, and let $\mathsf{m}=(m_{T,S})_{(T,S)\in\Pi(G)}$ be a sequence of integers indexed by $\Pi(G)$, constant on $G$-conjugacy classes of slices. Then the sequence $[\mathsf{m}]=(m_{T,S})_{(T,S)\in[\Pi(G)]}$ of representatives lies in the image of the ghost map~$\Phi$ if and only if, for any slice $(T,S)$ of $G$
$$\sum_{g\in N_G(T,S)/S}m_{\gen{gT},\gen{gS}}\equiv 0\;\big({\rm mod.}\,|N_G(T,S)/S|\big)\mpoint$$
\end{mth}
\pf Saying that $[\mathsf{m}]$ lies in the image of $\Phi$ is equivalent to saying that the element
$$s_\mathsf{m}=\sum_{(T,S)\in[\Pi(G)]}m_{T,S}\,\xi_{T,S}^G$$
of $\Q\Gringo(G)$ lies in $\Gringo(G)$, i.e. that it is a linear combination with integer coefficients of the elements $\gsect{V,U}_G$, for $(V,U)\in[\Pi(G)]$. Now by Theorem~\ref{slices and idempotents}
\begin{eqnarray*}
s_\mathsf{m}&=&\sum_{(T,S)\in \Pi(G)}\frac{|N_G(T,S)|}{|G|}m_{T,S}\,\xi_{T,S}^G\\
&=&\frac{1}{|G|}\sum_{(T,S)\in\Pi(G)}m_{T,S}\sum_{(V,U)\preceq (T,S)}|U|\mu_\Pi\big((V,U),(T,S)\big)\,\gsect{V,U}_G\\
&=&\frac{1}{|G|}\sum_{(V,U)\in\Pi(G)}|U|\Big(\sum_{(V,U)\preceq (T,S)}\mu_\Pi\big((V,U),(T,S)\big)m_{T,S}\Big)\gsect{V,U}_G\\
&=&\sum_{(V,U)\in[\Pi(G)]}\frac{1}{|N_G(V,U)/U|}\Big(\sum_{(V,U)\preceq (T,S)}\mu_\Pi\big((V,U),(T,S)\big)m_{T,S}\Big)\gsect{V,U}_G\mpoint
\end{eqnarray*}
Hence $s_\mathsf{m}\in\Im\,\Phi$ if and only if the number
$$\beta_{V,U}=\frac{1}{|N_G(V,U)/U|}\Big(\sum_{(V,U)\preceq (T,S)}\mu_\Pi\big((V,U),(T,S)\big)m_{T,S}\Big)$$
is an integer, for any slice $(V,U)$ of $G$.\par
With this notation, for any $(Y,X)\in \Pi(G)$
$$m_{Y,X}=\sumb{(V,U)\in\Pi(G)}{(Y,X)\preceq(V,U)}\beta_{V,U}|N_G(V,U)/U|\mpoint$$
Hence, setting $\sigma_{Y,X}=\sum_{g\in N_G(Y,X)/X}\limits m_{\gen{gY},\gen{gX}}$~:
\begin{eqnarray*}
\sigma_{Y,X}\!\!\!&\!=\!&\!\!\sum_{g\in N_G(Y,X)/X}\sumb{(V,U)\in\Pi(G)}{(\gen{gY},\gen{gX})\preceq(V,U)}\beta_{V,U}|N_G(V,U)/U|\\
\!\!&\!=\!&\!\!\sumc{g\in N_G(Y,X)}{(V,U)\in\Pi(G)}{(\gen{gY},\gen{gX})\preceq(V,U)}\frac{1}{|X|}\beta_{V,U}|N_G(V,U)/U|\\
\!\!&\!=\!&\!\!\sumc{(V,U)\in\Pi(G)}{(Y,X)\preceq (V,U)}{g\in U\cap N_G(Y,X)}\frac{1}{|X|}\beta_{V,U}|N_G(V,U)/U|\\
\end{eqnarray*}
\begin{eqnarray*}
\sigma_{Y,X}\!\!\!\!&\!=\!&\!\!\sumb{(V,U)\in [N_G(Y,X)\dom\Pi(G)]}{(Y,X)\preceq (V,U)}\frac{|N_G(Y,X)/X|}{|N_G(Y,X,V,U)|}\beta_{V,U}|N_G(V,U)/U||U\cap N_G(Y,X)|\\
\!\!&\!=\!&\!\!|N_G(Y,X)/X|\sumb{(V,U)\in [N_G(Y,X)\dom\Pi(G)]}{(Y,X)\preceq (V,U)}\beta_{V,U}\frac{|N_G(V,U)||U\cap N_G(Y,X)|}{|U||N_G(Y,X,V,U)|}\\
\!\!&\!=\!&\!\!|N_G(Y,X)/X|\sumb{(V,U)\in [N_G(Y,X)\dom\Pi(G)]}{(Y,X)\preceq (V,U)}\beta_{V,U}|N_G(V,U):U{\cdot} N_G(Y,X,V,U)|\mpoint
\end{eqnarray*}
Setting $\alpha_{Y,X}=\sigma_{Y,X}/|N_G(Y,X)/X|$, it follows that
$$\alpha_{Y,X}=\sumb{(V,U)\in N_G(Y,X)\dom\Pi(G)}{(Y,X)\preceq (V,U)}\beta_{V,U}|N_G(V,U):U{\cdot} N_G(Y,X,V,U)|\mpoint$$
Since $|N_G(V,U):U{\cdot} N_G(Y,X,V,U)|=1$ if $(Y,X)=(V,U)$, it follows that the transition matrix from the $\beta_{V,U}$'s to the $\alpha_{Y,X}$'s is triangular, with integer coefficients, and 1's on the diagonal. Hence it is invertible over $\Z$, and the numbers $\beta_{V,U}$ are all integers, for $(V,U)\in \Pi(G)$, if and only if the numbers $\alpha_{Y,X}$ are all integers, for $(Y,X)\in \Pi(G)$. This completes the proof.\findemo
\section{Prime spectrum}\label{Prime spectrum}
\begin{mth}{Notation} Let $p$ denote either 0 or a prime number.
\begin{itemize}
\item If $(T,S)\in\Pi(G)$, let  $I_{T,S,p}$ be the prime ideal of $\Gringo(G)$ defined as the kernel of the ring homomorphism 
$$\Gringo(G)\stackrel{\phi_{T,S}}{\longrightarrow}\Z\to \Z/p\Z\mvirg$$
where the right hand side map is the projection.
\item Let $\Theta(G)$ denote the set of triples $(T,S,p)$, where $(T,S)\in\Pi(G)$ is such that $|N_G(T,S)/S|\not\equiv 0 \;({\rm mod.}\;p)$.
\end{itemize}
\end{mth}
The group $G$ acts on $\Theta(G)$, by $^g(T,S,p)=({^gT},{^gS},p)$, for $g\in G$, and the ideal $I_{T,S,p}$ only depends on the $G$-orbit of $(T,S,p)$. Conversely~:
\begin{mth}{Proposition} \label{prime ideals}Let $I$ be a prime ideal of $\Gringo(G)$, and $R=\Gringo(G)/I$. Denote by $\phi:\Gringo(G)\to R$ the projection map, and denote by $p\geq 0$ the characteristic of $R$. Then $R\cong \Z/p\Z$ and~:
\begin{enumerate}
\item If $p=0$, there exists a slice $(T,S)$ of $G$ such that $\phi=\phi_{T,S}$, and $(T,S)$ is unique up to $G$-conjugation, with this property.
\item If $p>0$, there exists a slice $(T,S)$ of $G$ such that $\phi$ is the reduction modulo $p$ of $\phi_{T,S}$ and $N_G(T,S)/S$ is a $p'$-group, and $(T,S)$ is unique up to $G$-conjugation, with these properties.
\end{enumerate}
In particular, there exists a unique $(T,S,p)\in\Theta(G)$, up to conjugation, such that $I=I_{T,S,p}$.
\end{mth}
\pf Let $(T,S)$ be a slice of $G$, minimal for the relation $\preceq$, such that $\gsect{T,S}_G\notin I$. Then by Proposition~\ref{product formula}, for any $(Y,X)\in\Pi(G)$
\begin{eqnarray*}
\gsect{T,S}_G\gsect{Y,X}_G&=&\sum_{g\in[S\dom G/X]}\gsect{T\cap {^gY},S\cap{^gX}}_G\\
&\equiv&\sumc{g\in G/X}{S\leq{^gX}}{T\leq{^gY}}\gsect{T,S}_G \;\;({\rm mod.}\, I)\\
&=&\phi_{T,S}\big(\gsect{Y,X}_G\big)\gsect{T,S}_G\mpoint
\end{eqnarray*}
Since $I$ is prime, it follows that $\gsect{Y,X}_G-\phi_{T,S}(\gsect{Y,X}_G)1_{\Gringo(G)}\in I$. In particular $R=\Gringo(G)/I$ is generated by the image of $1_{\Gringo(G)}$, hence $R\cong \Z/p\Z$, where $p$ is the characteristic of $R$. Since $R$ is an integral domain, the number $p$ is either 0 or a prime.\begin{enumerate}
\item If $p=0$, then $R=\Z$, and $\phi=\phi_{T,S}$. And if $(T',S')\in\Pi(G)$ is such that $\phi_{T,S}=\phi_{T',S'}$, then both $\phi_{T,S}\big(\langle T',S'\rangle_G\big)$ and $\phi_{T',S'}\big(\langle T,S\rangle_G\big)$ are non zero. Then there exist elements $g,g'\in G$ such that $(T^g,S^g)\preceq (T',S')$ and $(T'^{g'},S'^{g'})\preceq (T,S)$, so $(T,S)$ and $(T',S')$ are conjugate in $G$.
\item If $p>0$, then $R=\Z/p\Z$, and $\phi$ is equal to the reduction of $\phi_{T,S}$ modulo~$p$. Since $\phi\big(\gsect{T,S}_G\big)=|N_G(T,S)/S|$ is non zero in $R$, it follows that $N_G(T,S)/S$ is a $p'$-group. If $(T',S')$ is another slice of $G$ such that $\phi$ is the reduction modulo $p$ of $\phi_{T',S'}$, and $N_G(T',S')/S'$ is a $p'$-group, then 
$$|N_G(T,S)/S|=\phi_{T,S}\big(\gsect{T,S}_G\big)\equiv \phi_{T',S'}\big(\gsect{T,S}_G\big)\;\;({\rm mod.}\,p)\mpoint$$ 
This is non zero. Similarly $|N_G(T,S)/S|\equiv \phi_{T,S}\big(\gsect{T',S'}_G\big)\;\;({\rm mod.}\,p)$ is non zero. In particular $\phi_{T',S'}\big(\gsect{T,S}_G\big)$ and $\phi_{T,S}\big(\gsect{T',S'}_G\big)$ are both non zero, and it follows as above that $(T,S)$ and $(T',S')$ are conjugate in~$G$.\findemo
\end{enumerate}
\begin{mth}{Notation} Let $p$ be a prime number. \begin{itemize}
\item Let $\Pi_p(G)$ denote the subset of $\Pi(G)$ consisting of the slices $(T,S)$ such that $N_G(T,S)/S$ is a $p'$-group.
\item For any $(T,S)\in\Pi(G)$, let $(T,S)\comp{p}$ denote the unique element $(V,U)$ of $\Pi_p(G)$, up to conjugation, such that $I_{T,S,p}=I_{V,U,p}$.
\end{itemize}
\end{mth}
\begin{mth}{Proposition} \label{closure} Let $p$ be a prime number. If $(T,S)$ is a slice of $G$, let $(T,S)\plus{p}$ denote a slice of the form $(PT,PS)$ of $G$, where $P$ is a Sylow $p$-subgroup of $N_G(T,S)$. \par
Define inductively an increasing sequence $(T_n,S_n)$ in $(\Pi(G),\preceq)$ by $(T_0,S_0)=(T,S)$, and $(T_{n+1},S_{n+1})=(T_n,S_n)\plus{p}$, for $n\in\N$. Then $(T,S)\comp{p}$ is conjugate to the largest term $(T_\infty,S_\infty)$ of the sequence $(T_n,S_n)$. 
\end{mth}
\pf The proof of Proposition~\ref{prime ideals} shows that $(T,S)\comp{p}$ is a minimal element $(V,U)$ of the poset $(\Pi(G),\preceq)$ such that
$$\phi_{T,S}(V,U)=|\{g\in G/U\mid (T^g,S^g)\preceq (V,U)\}|\not\equiv 0\;({\rm mod.} p)\mpoint$$
Thus one can assume that $(T,S)\preceq (V,U)$.  But $\phi_{T,S}\equiv \phi_{PT,PS}\;({\rm mod.}\;p)$ by Corollary~\ref{modp}, for any $p$-subgroup $P$ of $N_G(T,S)$, hence one can also assume that $(T,S)\plus{p}\preceq (V,U)$, and by induction, that $(T_\infty,S_\infty)\preceq (V,U)$. Moreover $\phi_{T_\infty,S_\infty}\equiv \phi_{V,U}\;({\rm mod.}\;p)$. As $N_G(T_\infty,S_\infty)/S_\infty$ is a $p'$-group, it follows that $(T_\infty,S_\infty)= (V,U)$, as was to be shown. \findemo
\begin{rem}{Remark} Let $(T,S)\in \Pi(G)$, and $(V,U)\in\Pi_p(G)$. It is easy to check, by induction on the integer $n$ such that $(T_n,S_n)=(T_\infty,S_\infty)$, that $(T,S)\comp{p}$ is equal to $(V,U)$ if and only if $T$ is a subnormal subgroup of $V$, if $S$ is a subnormal subgroup of $U$, if $|U:S|$ is a power of $p$, and if the set $T{\cdot}U$ is equal to $V$.
\end{rem}
\begin{mth}{Proposition} \label{ideal inclusion}Let $(T,S,p)$, $(T',S',p')$ be elements of $\Theta(G)$. Then $I_{T',S',p'}\subseteq I_{T,S,p}$ if and only if
\begin{itemize}
\item either $p'=p$ and the slices $(T',S')$ and $(T,S)$ are conjugate in~$G$.
\item or $p'=0$ and $p>0$, and the slices $(T',S')\comp{p}$ and $(T,S)$ are conjugate in~$G$.
\end{itemize}
\end{mth}
\pf Set $I=I_{T,S,p}$ and $I'=I_{T',S',p'}$. Then $\Gringo(G)/I'\cong \Z/p'\Z$ maps surjectively to $\Gringo(G)/I\cong\Z/p\Z$. Thus if $p=p'$, this projection map is an isomorphism, hence $I=I'$ and the slices $(T,S)$ and $(T',S')$ are conjugate in~$G$. And if $p\neq p'$, then $p'=0$ and $p>0$. The morphism $\phi_{T,S}$ is equal to the reduction modulo $p$ of the morphism $\phi_{T',S'}$. In other words $I_{T',S',p}=I_{T,S,p}$, hence $(T,S)$ is conjugate to $(T',S')\comp{p}$.\findemo
\begin{mth}{Corollary} \label{Zp components}Let $p$ be a prime number, and let $\Z_{(p)}$ be the localization of $\Z$ at the set $\Z-p\Z$. Let $\Theta_p(G)$ denote the subset of $\Theta(G)$ consisting of triples $(T,S,0)$, for $(T,S)\in\Pi(G)$, and $(T,S,p)$, for $(T,S)\in\Pi_p(G)$. Then~:
\begin{enumerate}
\item The prime ideals of the ring $\Z_{(p)}\Xi(G)$ are the ideals $\Z_{(p)}I_{T,S,q}$, for $(T,S,q)\in\Theta_p(G)$.
\item If $(T,S,q), (T',S',q')\in\Theta_p(G)$, then $\Z_{(p)}I_{T',S',q'}\subseteq \Z_{(p)}I_{T,S,q}$ if and only~if~:
\begin{itemize}
\item either $q=q'$, and the slices $(T,S)$ and $(T',S')$ are conjugate in~$G$.
\item or $q'=0$, $q=p$, and the slices $(T',S')\comp{p}$ and $(T,S)$ are conjugate in $G$.
\end{itemize}
\item The connected components of the spectrum of $\Z_{(p)}\Xi(G)$ are indexed by the conjugacy classes of $\Pi_p(G)$. The component indexed by $(T,S)\in\Pi_p(G)$ consists of a unique maximal element $\Z_{(p)}I_{T,S,p}$, and of the ideals $\Z_{(p)}I_{T',S',0}$, where $(T',S')\in\Pi(G)$ is such that $(T',S')\comp{p}$ is conjugate to $(T,S)$ in $G$.
\end{enumerate}
\end{mth} 
\pf The prime ideals of $\Z_{(p)}\Xi(G)$ are of the form $\Z_{(p)}I$, where $I$ is a prime ideal of $\Xi(G)$ such that $I\cap(\Z-p\Z)=\emptyset$. Equivalently $I=I_{T,S,0}$ or $I=I_{T,S,p}$. This proves Assertion~1. Now Assertion~2 follows from Proposition~\ref{ideal inclusion}, and Assertion~3 follows from Assertion~2.\findemo
\begin{mth}{Corollary} \label{pi idempotents}\begin{enumerate}
\item The primitive idempotents of the ring $\Z_{(p)}\Xi(G)$ are indexed by the conjugacy classes of $\Pi_p(G)$. The primitive idempotent $\eta_{V,U}^G$ indexed by $(V,U)\in\Pi_p(G)$ is equal to
$$\eta_{V,U}^G=\sumb{(T,S)\in[\Pi(G)]}{(T,S)\comp{p}=_G(T,S)}\xi_{T,S}^G\mpoint$$
\item Let $\pi$ be a set of prime numbers, and $\Z_{(\pi)}$ be the localization of $\Z$ relative to $\Z-\cup_{p\in\pi}p\Z$. Let $\mathcal{F}$ be a set of slices of $G$, invariant by $G$-conjugation, and $[\mathcal{F}]$ be a set of representatives of $G$-conjugacy classes of $\mathcal{F}$. Then the following conditions are equivalent~:
\begin{enumerate}
\item The idempotent 
$$\xi_\mathcal{F}^G=\sum_{(T,S)\in[\mathcal{F}]}\xi_{T,S}^G$$
of $\Q\Xi(G)$ lies in $\Z_{(\pi)}\Xi(G)$.
\item Let $(T,S)\in \Pi(G)$, and let $P$ be a $p$-subgroup of $N_G(T,S)$, for some $p\in\pi$. Then $(T,S)\in\mathcal{F}$ if and only if $(PT,PS)\in\mathcal{F}$. 
\end{enumerate}
\end{enumerate}
\end{mth}
\pf Let $\mathcal{F}$ be a set of slices of $G$, invariant by $G$-conjugation, and $[\mathcal{F}]$ be a set of representatives of $G$-conjugacy classes of $\mathcal{F}$. The idempotent
$$\xi_\mathcal{F}^G=\sum_{(T,S)\in[\mathcal{F}]}\xi_{T,S}^G$$
of $\Q\Xi(G)$ lies in $\Z_{(p)}\Xi(G)$, for some prime $p$, if and only if there exists an integer $m$, not divisible by $p$, such that $u=m\xi_\mathcal{F}^G\in\Xi(G)$. Let $(T,S)\in\Pi(G)$, and let $P$ be a $p$-subgroup of $N_G(T,S)$. The integer $\phi_{T,S}(u)$ is equal to $m$ if $(T,S)\in\mathcal{F}$, and to 0 otherwise. Hence it is coprime to $p$ if and only if $(T,S)\in\mathcal{F}$. Since $\phi_{T,S}$ and $\phi_{PT,PS}$ are congruent modulo $p$, it follows that $(T,S)\in\mathcal{F}$ if and only if $(PT,PS)\in\mathcal{F}$. \par
Hence if $(T,S)$ and $(T',S')$ are slices of $G$ such that $(T,S)\comp{p}=_G(T',S')\comp{p}$, then $(T,S)\in\mathcal{F}$ if and only if $(T',S')\in\mathcal{F}$. Thus $\mathcal{F}$ is a disjoint union of sets of the form 
$$E_{V,U}=\{(T,S)\in\Pi(G)\mid (T,S)\comp{p}=_G(V,U)\}\mvirg$$
for some slices $(V,U)\in\Pi_p(G)$. In other words the idempotent $\xi_\mathcal{F}^G$ is a sum of some idempotents $\eta_{V,U}^G$, for $(V,U)\in\Pi_p(G)$. \par
But the primitive idempotents of the ring $\Z_{(p)}\Xi(G)$ are in one to one correspondence with the connected components of its spectrum, which precisely are indexed by the conjugacy classes of $\Pi_p(G)$. It follows that $\xi_\mathcal{F}^G=\eta_{V,U}^G$ is equal to the idempotent corresponding to the component indexed by $(V,U)$, for any $(V,U)\in\Pi_p(G)$. This proves Assertion~1. This also proves Assertion~2 in the case where $\pi$ consists of a single prime number.\par
For the general case, observe that $\xi_{\mathcal{F}}^G$ lies in $\Z_{(\pi)}\Xi(G)$ if and only if it lies in $\Z_{(p)}\Xi(G)$, for any $p\in \pi$.\findemo
\begin{mth}{Theorem} \label{connected Xi}Let $G$ be a finite group. 
\begin{enumerate}
\item Let $\sim$ denote the finest equivalence relation on the set $\Pi(G)$ such that for any $(T,S), (T',S')\in\Pi(G)$,
$$\exists p,\;(T,S)\comp{p}=_G(T',S')\comp{p}\implies (T,S)\sim (T',S')\mpoint$$
Then the primitive idempotents of $\Xi(G)$ are indexed by the equivalence classes of $\Pi(G)$ for the relation $\sim$. The idempotent $\xi_\mathcal{C}^G$ indexed by the component $\mathcal{C}$ is equal to
$$\xi_\mathcal{C}^G=\sum_{(T,S)\in[\mathcal{C}]}\xi_{T,S}^G\mvirg$$
where $[\mathcal{C}]$ is a set of representatives of $G$-conjugacy classes in $\mathcal{C}$.
\item The prime spectrum of $\Xi(G)$ is connected if and only if $G$ is solvable.
\end{enumerate}
\end{mth}
\pf Assertion 1 follows from Corollary~\ref{pi idempotents}, applied to the set $\pi$ of all primes.\par
For Assertion 2, observe that the spectrum of $\Xi(G)$ is connected if and only if 1 is a primitive idempotent of $\Xi(G)$. This means that for any two slices $(T,S)$ and $(T',S')$ of $G$, there exists a sequence $(T_i,S_i)$ of slices, for $i\in\{0,\ldots,n\}$, and a sequence $p_i$ of prime numbers, for $i\in\{0,\ldots,n-1\}$, such that
$$(T,S)=(T_0,S_0)\stackrel{p_0}{\sim}(T_1,S_1)\stackrel{p_1}{\sim}\ldots\stackrel{p_{n-2}}{\sim}(T_{n-1},S_{n-1})\stackrel{p_{n-1}}{\sim}(T_n,S_n)=(T',S')\mvirg$$
where the notation $(Y,X)\stackrel{p}{\sim}(Y',X')$ means that $(Y,X)\comp{p}=_G(Y',X')\comp{p}$. \par
But clearly, if $(Y,X)\stackrel{p}{\sim}(Y',X')$, then the slices $\big(O^p(Y),O^p(X)\big)$ and $\big(O^p(Y'),O^p(X')\big)$ are conjugate in $G$. Hence the slices $\big(D^\infty(Y),D^\infty(X)\big)$ and $\big(D^\infty(Y'),D^\infty(X')\big)$ are conjugate in $G$, where $D^\infty(H)$ denotes the last term in the derived series of the group $H$. \par
If the spectrum of $G$ is connected, taking $(T,S)=(G,G)$ and $(T',S')=(\un,\un)$, it follows that $D^\infty(G)=\un$, i.e. that $G$ is solvable.\par
Conversely, if $G$ is solvable, let $(T,S)$ be a slice of $G$. Then $S$ is solvable, and there exists a prime $p$ such that $O^p(S)<S$. Let $P$ be a Sylow $p$-subgroup of $S$. Then $P\leq N_G(T,S)$, and $(T,S)=\big(PT,PO^p(S)\big)$. Thus $(T,S)\stackrel{p}{\sim}\big(T,O^p(S)\big)$. By induction, there is a sequence of prime numbers $p_i$, for $i\in\{1,\ldots,k\}$, such that
$$(T,S)=(T,S_0)\stackrel{p_0}{\sim}(T,S_1)\stackrel{p_1}{\sim}\ldots\stackrel{p_{k-1}}{\sim}(T,S_{k})\stackrel{p_{k}}{\sim}(T,\un)\mvirg$$
where $S_{i+1}=O^{p_i}(S_i)$.\par
Now $T$ is solvable, so there exists a prime $q$ such that $O^q(T)<T$. If $Q$ is a Sylow $q$-subgroup of $T$, then $Q\leq N_G(T,\un)$, and $(T,\un)\stackrel{q}{\sim}\big(O^q(T),\un\big)$. Hence there exists a sequence of primes $q_j$, for $j\in\{0,\ldots,l\}$, such that
$$(T,\un)=(T_0,\un)\stackrel{q_0}{\sim}(T_1,\un)\stackrel{q_1}{\sim}\ldots\stackrel{q_{l-1}}{\sim}(T_{l},\un)\stackrel{q_{l}}{\sim}(\un,\un)\mvirg$$
where $T_{j+1}=O^{q_j}(T_j)$. This shows that the spectrum of $\Xi(G)$ is connected.\findemo
\section{Unit group}\label{units}
\npar When $G$ is a finite group, denote by $\Xi(G)^\times$ the group of invertible elements of the ring $\Xi(G)$. It follows from Theorem~\ref{Gringo basis} that the restricted ghost map yields an injective group homomorphism 
$$\Phi^\times:\Xi(G)^\times\hookrightarrow \prod_{(T,S)\in\Pi(G)}\Z^\times\mpoint$$
The following lemma is a straightforward consequence of the existence of this injective group homomorphism~:
\begin{mth}{Lemma} \label{units equivalence}Let $G$ be a finite group, and let $u\in\Xi(G)$. The following conditions are equivalent~:
\begin{enumerate}
\item $u\in \Xi(G)^\times$.
\item $\phi_{T,S}(u)\in\{\pm 1\}$, for any $(T,S)\in\Pi(G)$.
\item $u^2=1$.
\end{enumerate}
In particular $\Xi(G)^\times$ is a finite elementary abelian 2-group.\par
\end{mth}
The main motivation in considering the group $\Xi(G)^\times$ lies in the following proposition, which extends a theorem of tom Dieck (\cite{tomdieckgroups} Proposition 1.5.1) about the unit group of the usual Burnside ring~:
\begin{mth}{Proposition} \label{Feit-Thompson}Feit-Thompson's theorem is equivalent to the statement that, if $G$ has odd order, then $\Xi(G)^\times=\{\pm 1\}$.
\end{mth}
\pf The first observation is that for any finite group $G$, by Lemma~\ref{units equivalence}, the correspondences $u\mapsto \frac{1-u}{2}$ and $e\mapsto 1-2e$ are mutually inverse bijections between $\Xi(G)^\times$ and the set of idempotents $e\in\Q \Xi(G)$ such that $2e\in\Xi(G)$.\par
Now  Theorem~\ref{slices and idempotents} shows that $|G|e\in\Xi(G)$, for any idempotent $e$ of $\Q \Xi(G)$. Hence if $G$ has odd order, and if $e$ is an idempotent of $\Xi(G)$ such that $2e\in \Xi(G)$, then $(2,|G|)e=e\in\Xi(G)$. Thus if $|G|$ is odd, the set $\Xi(G)^\times$ is in one to one correspondence with the set of idempotents of the ring $\Xi(G)$. By a standard argument from commutative ring theory, this set is in one to one correspondence with the set of connected components of the spectrum of $\Xi(G)$. It follows that if $G$ has odd order, the spectrum of $G$ is connected if and only if $\Xi(G)^\times=\{\pm 1\}$. By Theorem~\ref{connected Xi}, this is equivalent to saying that $G$ is solvable.\findemo
The following theorem is an analogue of Yoshida's characterization (\cite{yoshidaunit}) of the unit group of the usual Burnside ring~:
\begin{mth}{Theorem} \label{units characterization}Let $G$ be a finite group, and let $\mathsf{m}=(m_{T,S})_{(T,S)\in\Pi(G)}$ be a sequence of integers in $\{\pm 1\}$ indexed by $\Pi(G)$, constant on $G$-conjugacy classes of slices. Then the sequence $[m]= (m_{T,S})_{(T,S)\in[\Pi(G)]}$ of representatives lies in the image of the restricted ghost map $\Phi^\times$ if and only if for any $(T,S)\in\Pi(G)$, the map
$$g\in N_G(T,S)/S\mapsto m_{\gen{gT},\gen{gS}}/m_{T,S}\in\{\pm 1\}$$
is a group homomorphism.
\end{mth}
\pf Let $X\stackrel{f}{\to}Y$ be a morphism of finite $G$-sets. It follows from Lemma~\ref{phiTS} that for any $(T,S)\in \Pi(G)$, the monoid of endomorphisms of $G/S\to G/T$ in the category $\mor{G}$ is actually a group, isomorphic to $N_G(T,S)/S$. This group acts on the set of morphisms from $G/S\to G/T$ to $X\stackrel{f}{\to}Y$, by pre-composition~: if
$$\xymatrix{
G/S\ar[r]\ar[d]^-\alpha&G/T\ar[d]^-{\beta}\\
X\ar[r]^-f&Y
}
$$
is a morphism in $\mor{G}$, and if $g\in N_G(T,S)$, the morphism $^g(\alpha,\beta)=({^g\alpha},{^g\beta})$ is defined by $(^g\alpha)(xS)=xgS$ and $(^g\beta)(xT)=xgT$, for any $g\in T$. The morphism $(\alpha,\beta)$ is invariant under $g\in N_G(T,S)$ if and only if $\alpha(gS)=\alpha(S)$, i.e. if $(\alpha,\beta)$ factors as
$$\xymatrix{
G/S\ar[r]\ar[d]&G/gT\ar[d]\\
G/\gen{gS}\ar[r]\ar[d]^-{\sur{\alpha}}&G/\gen{gT}\ar[d]^-{\sur{\beta}}\\
X\ar[r]^-f&Y
}
$$
It follows that the number of fixed points of $g\in N_G(T,S)/S$ on the set of homomorphisms from $G/S\to G/T$ to $X\stackrel{f}{\to}Y$, i.e. the value at $gS$ of the corresponding permutation character $\theta_{T,S}$ of $N_G(T,S)/S$,  is equal to $\phi_{\gen{gT},\gen{gS}}(X\stackrel{f}{\to}Y)$. This shows more generally that for any $u\in \Xi(G)$, the correspondence 
$$\theta_{T,S}:g\in N_G(T,S)/S\mapsto \phi_{\gen{gT},\gen{gS}}(u)$$
is a generalized character of the group $N_G(T,S)/S$.\par
Now if $u\in\Xi(G)$, this generalized character has all its values in $\{\pm 1\}$. It follows that $\gsect{\theta_{T,S},\theta_{T,S}}_G=1$, hence $\theta_{T,S}$ is up to a sign equal to an irreducible character of $N_G(T,S)/S$, of degree 1. Hence $\theta_{T,S}/\theta_{T,S}(1)$ is a group homomorphism from $N_G(T,S)/S$ to $\{\pm 1\}$. Thus if $\mathsf{m}\in \Im(\Phi^\times)$, then the map $g\in N_G(T,S)/S\mapsto m_{\gen{gT},\gen{gS}}/m_{T,S}\in\{\pm 1\}$ is a group homomorphism.\par
Conversely, if $\mathsf{m}=(m_{T,S})_{(T,S)\in\Pi(G)}$ is a $G$-invariant sequence with values in $\{\pm 1\}$, such that for any $(T,S)\in \Pi(G)$, the map $g\in N_G(T,S)/S\mapsto m_{\gen{gT},\gen{gS}}/m_{T,S}\in\{\pm 1\}$ is a group homomorphism $\theta_{T,S}$, then 
$$\sum_{g\in N_G(T,S)/S}m_{\gen{gT},\gen{gS}}=m_{T,S}|N_G(T,S)/S|\,\gsect{\theta_{T,S},1}_G$$
is an integer multiple of $|N_G(T,S)/S|$. By Theorem~\ref{ghost image}, it follows that $[\mathsf{m}]=\Phi(u)$, for some $u\in\Xi(G)$. Then $u\in\Xi(G)^\times$, by Lemma~\ref{units equivalence}, and this completes the proof.\findemo
\npar In the case of the usual Burnside ring $B(G)$, the correspondence sending a finite group $G$ to the unit group $B(G)^\times$ can be endowed with a structure of biset functor, using {\em tensor induction} (see e.g. \cite{burnsideunits}, Proposition 5.5). One may ask whether a similar structure exists for the group of units of the section Burnside ring~: 
\begin{mth}{Proposition} \label{not biset functor}The correspondence sending a finite group $G$ to $\Xi(G)^\times$ cannot be endowed with a structure of biset functor.
\end{mth}
\pf If such a structure exists, one may view it as a biset functor $\Xi^\times$ with values in $\F_2$-vector spaces. Obviously $\Xi^\times(\un)=\F_2$. This shows that there are two subfunctors $F_2\subset F_1$ of $\Xi^\times$ such that $F_1/F_2$ is isomorphic to the simple functor $S_{\un,\F_2}$. \par
An easy computation (a special case of Theorem~\ref{unit group elemab} in the Appendix) shows moreover that $\Xi^\times(C_2)\cong (\F_2)^3$, where $C_2$ is a group of order 2. As $S_{\un,\F_2}(C_2)\cong (\F_2)^2$ (see e.g. \cite{bisetfunctors} Proposition~4.4.6), this shows that either $(F/F_1)(C_2)\cong \F_2$ and $F_2(C_2)=\zero$, or $F(C_2)=F_1(C_2)$ and $F_2(C_2)\cong \F_2$. In the first case, there are subfunctors $F_3$ and $F_4$ of $F$ with $F_1\subseteq F_4\subset F_3\subseteq F$ and $F_3/F_4\cong S_{C_2,\F_2}$. In the latter case, there are subfunctors $F_4\subset F_3$ of $F_2$ such that $F_3/F_4\cong S_{C_2,\F_2}$. In any case
$$\dim_{\F_2}\Xi^\times\big((C_2)^2\big)\geq \dim_{\F_2}S_{\un,\F_2}\big((C_2)^2\big)+\dim_{\F_2}S_{C_2,\F_2}\big((C_2)^2\big)\mpoint$$
Now $\dim_{\F_2}S_{\un,\F_2}\big((C_2)^2\big)=4$ by \cite{bisetfunctors} Proposition~4.4.6 or Corollary~10.5.6, and $\dim_{\F_2}S_{C_2,\F_2}\big((C_2)^2\big)=5$ by inspection, using \cite{bisetfunctors} Proposition~4.4.6. Thus $\dim_{\F_2}\Xi^\times\big((C_2)^2\big)\geq 9$.\par
When $G$ is a finite abelian group, the unit group $\Xi(G)^\times$ can be determined explicitly (see Theorem~\ref{unit group elemab}). In particular
\begin{equation}\label{klein}\Xi\big((C_2)^2\big)^\times\cong(\F_2)^7\mpoint
\end{equation}
This contradiction shows that the biset functor $\Xi^\times$ cannot exist.\findemo
\begin{rem}{Remark} Boltje and Pfeiffer (\cite{boltje-pfeiffer}) have described an efficient algorithm to compute the unit group of the ordinary Burnside ring. A straightforward adaptation of this algorithm to the ring $\Xi(G)$ allows for a quick computation of the group $\Xi(G)^\times$, for not too large finite groups $G$, using GAP4 (\cite{GAP4}). These computations agree in particular with~\ref{klein}.
\end{rem}
\vspace{4ex}\par
\pagebreak[4]
\centerline{\Large\bf \romain{2} - The section Burnside ring}
\section{Galois morphisms of $G$-sets}\label{Galois morphisms}
\begin{mth}{Definition} Let $G$ be a group. A morphism $f:X\to Y$ of $G$-sets is {\em a Galois morphism} if for any $x,x'\in X$ such that $f(x)=f(x')$, there exists $\varphi\in\Aut_G(X)$ such that $f\circ\varphi=f$ and $\varphi(x)=x'$.
\end{mth} 
\begin{rem}{Example} \label{injective => galois}Any injective morphism of $G$-sets is a Galois morphism, for trivial reasons. 
\end{rem}
\begin{mth}{Proposition} \label{stabilizers} Let $f:X\to Y$ be a morphism of $G$-sets. The following conditions are equivalent:
\begin{enumerate}
\item $f$ is a Galois morphism.
\item $\forall x,x'\in X$, $f(x)=f(x')\implies G_x=G_{x'}$.
\end{enumerate}
\end{mth}
\pf Suppose first that $f$ is a Galois morphism. Then if $f(x)=f(x')$, there exists a $G$-automorphism $\varphi$ of $X$ such that $\varphi(x)=x'$. This implies $G_x\leq G_{x'}$, hence $G_x=G_{x'}$ by symmetry. Thus Condition 1 implies Condition~2.\par
Conversely, suppose that Condition 2 holds, and let $x,x'\in X$ with $f(x)=f(x')$. There are two cases~:\begin{itemize}
\item Either $x$ and $x'$ are in the same $G$-orbit $\omega$. Let $Z=X-\omega$. Define $\varphi:X\to X$ by $\varphi(t)=t$ if $t\in Z$ and $\varphi(ux)=ux'$ if $u\in G$. This is a well defined $G$-automorphism of $X$, since $G_x=G_{x'}$. Clearly $\varphi(x)=x'$, and $f\circ\varphi(t)=f(t)$ if $t\in Z$. Moreover 
$$f\circ\varphi(ux)=f(ux')=uf(x')=uf(x)=f(ux)$$
for any $u\in G$. Hence $f\circ\varphi=f$.
\item If $x$ and $x'$ are in different $G$-orbits $\omega$ and $\omega'$, respectively, let $Z=X-(\omega\sqcup\omega')$. Define $\varphi:X\to X$ by $\varphi(t)=t$ if $t\in Z$, and $\varphi(ux)=ux'$ and $\varphi(ux')=ux$, for any $u\in G$. This is a well defined $G$-automorphism of~$X$, since $G_x=G_{x'}$. Clearly $\varphi(x)=x'$, and $f\circ\varphi(t)=f(t)$ if $t\in Z$. Moreover 
$$f\circ\varphi(ux)=f(ux')=uf(x')=uf(x)=f(ux)$$ 
for any $u\in G$. Also 
$$f\circ\varphi(ux')=f(ux)=uf(x)=uf(x')=f(ux')\mvirg$$ 
for any $u\in G$. Hence $f\circ\varphi=f$.
\end{itemize}
It follows that $f$ is a Galois morphism, and Condition 2 implies Condition~1.\findemo
\begin{mth}{Corollary} \label{criterion}Let $f:X\to Y$ be a morphism of $G$-sets. Then $f$ is a Galois morphism if and only if for any $y\in f(X)$, there exists a normal subgroup $N_y$ of $G_y$ such that $G_x=N_y$ for any $x\in f^{-1}(y)$.
\end{mth}
\pf Suppose first that $f$ is a Galois morphism. If $x\in X$, set $y=f(x)$, and choose $g\in G_y$. Then $f(gx)=gy=y=f(x)$, thus $G_{gx}={^gG_x}=G_x$, thus $G_x\normal G_y$. Moreover $G_x$ does not depend on $x\in f^{-1}(y)$, by Proposition~\ref{stabilizers}.\par
Conversely, suppose that for any $y\in f(X)$, there exists $N_y\normal G_y$ such that $G_x=N_y$ for any $x\in f^{-1}(y)$. Then obviously $G_x=G_{x'}=N_y$ if $f(x)=f(x')=y$, so $f$ is a Galois morphism, by Proposition~\ref{stabilizers}.\findemo
\begin{rem}{Remark} In particular, when $(T,S)$ is a slice of $G$, the projection morphism $G/S\to G/T$ is a Galois morphism of $G$-sets if and only if $S\normal T$, i.e. if $(T,S)$ is a {\em section} of $G$.
\end{rem}

\begin{mth}{Lemma} \label{injection} Let $f:X\to Y$ be a morphism of $G$-sets, and $j:Y\hookrightarrow Z$ be an injective morphism of $G$-sets. Then $f$ is a Galois morphism if and only if $j\circ f$ is a Galois morphism.
\end{mth}
\pf This is straightforward.\findemo
\begin{mth}{Lemma}\label{base change} Let 
$$\xymatrix{X\ar[r]^-a\ar[d]_-b&Y\ar[d]^-c\\
Z\ar[r]_-d&T
}
$$
be a cartesian square of $G$-sets. If $c$ is a Galois morphism, then $b$ is a Galois morphism.
\end{mth}
\pf Suppose that $x,x'\in X$ are such that $b(x)=b(x')$. Then 
$$ca(x)=db(x)=db(x')=ca(x')\mpoint$$
If $c$ is a Galois morphism, it follows that $G_{a(x)}=G_{a(x')}$. Let $g\in G_x$. Then $g\in G_{a(x)}=G_{a(x')}$, and $g\in G_{b(x)}=G_{b(x')}$, thus
\begin{eqnarray*}
b(x')&=&gb(x')=b(gx')\\
a(x')&=&ga(x')=a(gx')\mpoint\\
\end{eqnarray*}
It follows that $x'=gx'$, so $G_x\leq G_{x'}$, and $G_x=G_{x'}$ by symmetry. Hence $f$ is a Galois morphism, by Proposition~\ref{stabilizers}. \findemo
\begin{rem}{Remark} In particular, any morphism isomorphic to a Galois morphism (that is, when $a$ and $d$ are isomorphisms) is a Galois morphism.
\end{rem}
\begin{mth}{Corollary} \label{restriction to subdomain} Let $f:X\to Y$ be a Galois morphism of $G$-sets. If $X_1$ is a $G$-subset of $X$, the restricted morphism $f_{|X_1}:X_1\to f(X_1)$ is a Galois morphism of $G$-sets.
\end{mth}
\pf This follows from Lemma~\ref{base change} and Lemma~\ref{injection}, since the square
$$\xymatrix{X_1\ar[r]^-{i}\ar[d]_-{f_1}&X\ar[d]^-f\\
f(X_1)\ar[r]_-j&Y
}
$$
is cartesian, where $i$ and $j$ are the injection maps, and $f_1=f_{\mid X_1}$.\findemo

\begin{mth}{Proposition}\label{base} A morphism $f:X\to Y_1\sqcup Y_2$ is a Galois morphism if and only if the restricted morphisms $f^{-1}(Y_1)\to Y_1$ and  $f^{-1}(Y_2)\to Y_2$ are Galois morphisms.
\end{mth}
\pf Set $X_i=f^{-1}(Y_i)$, and denote by $f_i:X_i\to Y_i$ the restriction of $f$, for $i=1,2$. Assume first that $f$ is a Galois morphism. If $x,x'\in X_1$ have the same image under $f_1$, then they have the same image under $f$, and $G_x=G_{x'}$. So $f_1$ is a Galois morphism, and $f_2$ is also a Galois morphism, by symmetry.\par
Conversely, suppose that $f_1$ and $f_2$ are Galois morphisms. If $x,x'\in X$ are such that $f(x)=f(x')\in Y_1$, then $x,x'\in X_1$, and $f_1(x)=f_1(x')$. Thus $G_x=G_{x'}$, as $f_1$ is a Galois morphism. If $f(x)=f(x')\in Y_2$, the argument is similar, with $f_1$ replaced by $f_2$. In any case $G_x=G_{x'}$, and $f$ is a Galois morphism, by Proposition~\ref{stabilizers}.\findemo
\begin{mth}{Proposition} \label{biset} Let $G$ and $H$ be finite groups, and $U$ be an $(H,G)$-biset. If $f:X\to Y$ is a Galois morphism of $G$-sets, then 
$$U\times_Gf:U\times_GX\to U\times_GY\;\;\hbox{and}\;\;U\circ_Gf:U\circ_GX\to U\circ_GY$$
are Galois morphisms of $H$-sets.
\end{mth}
\pf For $U\times_Gf$, let $(u,_{_G}x)$ and $(u',_{_G}x')$ be elements of $U\times_GX$ having the same image under $U\times_Gf$. This means that $\big(u,_{_G}f(x)\big)=\big(u',_{_G}f(x')\big)$, i.e. that there exist $g\in G$ such that $ug^{-1}=u'$ and $gf(x)=x'$. Thus $f(gx)=f(x')$, and $G_{gx}=G_{x'}$, since $f$ is a Galois morphism. \par
Now if $h\in H$ stabilizes $(u,_{_G}x)$, i.e. if $(hu,_{_G}x)=(u,_{_G}x)$, there exists $a\in G$ such that $hu=ua$ and $a^{-1}x=x$. Hence
$$h(u',_{_G}x')=(hu',_{_G}x')=(hug^{-1},_{_G}x')=(uag^{-1},_{_G}x')=(u,_{_G}ag^{-1}x')\mpoint$$
Since $a\in G_x$, it follows that $^ga\in G_{gx}=G_{x'}$, so $ag^{-1}x'=g^{-1}x'$, and 
$$h(u',_{_G}x')=(u,_{_G}g^{-1}x')=(ug^{-1},_{_G}x')=(u',_{_G}x')\mpoint$$
This shows that $h$ stabilizes $(u',_{_G}x')$. By symmetry, the stabilizers of $(u,_{_G}x)$ and $(u',_{_G}x')$ in $H$ are equal, and $U\times_Gf$ is a Galois morphism of $H$-sets.\par
Now recall from \cite{doublact2} that $U\circ_GX$ is the $H$-subset of $U\times_GX$ defined by
$$U\circ_GX=\{(u,_{_G}x)\in U\times_GX\mid\forall g\in G,\;ug=g\implies gx=x\}\mvirg$$
and that the map $U\circ_Gf$ is the restriction of $U\times_Gf$ to $U\circ_GX$. By Lemma~\ref{injection} and Corollary~\ref{restriction to subdomain}, the morphism $U\circ_Gf$ is also a Galois morphism of $H$-sets.\findemo
\begin{mth}{Corollary} \label{restriction} Let $f:X\to Y$ be a Galois morphism of $G$-sets. If $H$ is a subgroup of $G$, the restriction $\Res_H^Gf:\Res_H^GX\to\Res_H^GY$ is a Galois morphism of $H$-sets.
\end{mth}
\pf  In Proposition~\ref{biset}, set $U=G$, viewed as and $(H,G)$-biset for left and right multiplication.\findemo
Let $G$ and $H$ be groups. A morphism $f:U\to U'$ of $(H,G)$-bisets is called a Galois morphism if it is a Galois morphism of $(H\times G\op)$-sets. Then~:
\begin{mth}{Proposition} Let $G$, $H$ and $K$ be groups. If $f:U\to U'$ is a Galois morphism of $(H,G)$-bisets, and $g:V\to V'$ is a Galois morphism of $(K,H)$-bisets, then $g\times_Hf:V\times_HU\to V'\times_HU'$ is a Galois morphism of $(K,G)$-bisets.
\end{mth}
\pf Let $(v,_{_H}u)$ and $(v',_{_H}u')$ be elements of $V\times_HU$ with the same image under $g\times_Hf$. It means that there exists $h\in H$ such that
$$g(v')=g(v)h=g(vh)\;\;\hbox{and}\;\;f(u')=h^{-1}f(u)=f(h^{-1}u)\mpoint$$
As $g$ is a Galois morphism, the stabilizers of $v'$ and $vh$ in $K\times H\op$ are equal. Similarly, as $f$ is a Galois morphism, the stabilizers of $u'$ and $h^{-1}u$ in $H\times G\op$ are equal.\par
Now let $(k,g)\in K\times G\op$ such that $k(v',_{_H}u')g^{-1}=(v',_{_H}u')$. It means that there exists $a\in H$ such that $kv'a^{-1}=v'$ and $au'g^{-1}=u'$. Hence
$$kvha^{-1}=vh\;\;\hbox{and}\;\;ah^{-1}ug^{-1}=h^{-1}u\mpoint$$
It follows that
\begin{eqnarray*}
k(v,_{_H}u)g^{-1}&=&(kv,_{_H}ug^{-1})=(v{hah^{-1}},_{_H}ug^{-1})\\
&=&(v,_{_H}{hah^{-1}}ug^{-1})=(v,_{_H}hh^{-1}u)=(v,_{_H}u)\mpoint
\end{eqnarray*}
By symmetry, the stabilizers of $(v,_{_H}u)$ and $(v',_{_H}u')$ in $K\times G\op$ are equal. Thus $(g\times_Hf)$ is a Galois morphism of $(K,G)$-bisets.\findemo
\begin{mth}{Corollary} \label{product} Let $G$ and $H$ be groups. If $f:X\to Y$ is a Galois morphism of $G$ sets and $g:Z\to T$ is a Galois morphism of $H$-sets, then $f\times g: X\times Z\to Y\times T$ is a Galois morphism of $(G\times H)$-sets.
\end{mth}
\pf Consider $f$ as a morphism of $(G,\un)$-bisets, and $g$ as a morphism of $(\un,H\op)$-bisets.\findemo
\begin{mth}{Notation} 
Let $\galmor{G}$ denote the full subcategory of $\mor{G}$ consisting of Galois morphisms of $G$-sets.
\end{mth}
\begin{mth}{Notation} Let $X\stackrel{f}{\to} Y$ be a morphism of $G$-sets. For $x\in X$, set
$$G_x^{f}=\gen{G_{z}\mid \forall z\in X,\;f(z)=f(x)}\mpoint$$
Let $\sim_f$ be the relation on $X$ defined by
$$x\sim_f x'\Leftrightarrow \exists g\in G_x^{f},\;gx=x'\mpoint$$
\end{mth}
\begin{mth}{Lemma} With this notation~:
\begin{enumerate}
\item The relation $\sim_f$ is an equivalence relation on $X$. Let $X_f\Gal$ denote the set of equivalence classes, and let $\gamma_{X,f}:X\to X\Gal_f$ denote the projection map.
\item If $x,x'\in X$ and $x\sim_f x'$, then $gx\sim_fgx'$ for any $g\in G$. Hence there exists a unique structure of $G$-set on $X\Gal_f$ such that $\gamma_{X,f}$ is a morphism of $G$-sets.
\item There is a unique map $f\Gal:X_f\Gal\to Y$ such that the diagram
$$\xymatrix{
X\ar[r]^-f\ar[d]_{\gamma_{X,f}}&Y\\
X_f\Gal\ar[ur]_-{f\Gal}
}
$$
is commutative.
\end{enumerate}
\end{mth} 
\pf For Assertion 1, the relation $\sim_f$ is clearly reflexive, since $G_x\leq G_x^f$. Now if $x,x'\in X$ and $x\sim_fx'$, let $g\in G_x^f$ such that $gx=x'$. There exist $r\in\N$, elements $z_1,\ldots, z_r$ in $X$ and elements $g_1,\ldots,g_r$ of $G$ such that $g_iz_i=z_i$ and $f(z_i)=f(x)$, for $i=1,\ldots,r$, and such that $g=g_1\cdots g_r$. It follows that
\begin{eqnarray*}
f(x')&=&g_1\cdots g_rf(x)=g_1\cdots g_rf(z_r)=g_1\cdots g_{r-1}f(g_rz_r)=g_1\cdots g_{r-1}f(z_r)\\
&=&g_1\cdots g_{r-1}f(x)\mvirg
\end{eqnarray*}
thus $f(x')=f(x)$, by induction on $r$. It follows that $G_{x'}^f=G_x^f$, hence that $x'\sim_fx$, since $x=g^{-1}x'$. So the relation $\sim_f$ is symmetric.\par
Finally $\sim_f$ is transitive~: if $x,x',x''\in X$, if $x\sim_fx'$ and $x'\sim_fx''$, there exist $g\in G_x^f$ and $g'\in G_{x'}^f$ such that $gx=x'$ and $g'x'=x''$. But the previous argument shows that $G_x^f=G_{x'}^f=G_{x''}^f$. Hence $g'g\in G_x^f$, and $x\sim_fx''$, since $g'gx=x''$.\par
Assertion 2 follows from the straightforward fact that $^g(G_x^f)=G_{gx}^f$ for any $x\in X$ and any $g\in G$. Assertion~3 follows as well, since $x\sim_fx'$ implies $f(x)=f(x')$.\findemo
\begin{mth}{Proposition} \label{adjunction}Let $G$ be a finite group.
\begin{enumerate}
\item Let $X\stackrel{f}{\to}Y$ be a morphism of $G$-sets. Then the morphism $X\Gal_f\stackrel{f\Gal}{\longrightarrow}Y$ is a Galois morphism of $G$-sets.
\item If 
$$\xymatrix{
X\ar[r]^-f\ar[d]_-\alpha&Y\ar[d]^\beta\\
A\ar[r]^-a&B
}
$$
is a morphism in the category $\mor{G}$, and if $A\stackrel{a}{\to}B$ is a Galois morphism of $G$-sets, then there exists a unique morphism of $G$-sets $\widetilde{\alpha}:X_f\Gal\to A$ such that the diagram
$$\xymatrix{
X\ar[r]^-f\ar[d]^-{\gamma_{X,f}}\ar@/_5ex/[dd]_-\alpha&Y\ar@{=}[d]\\
X\Gal_f\ar[d]^-{\widetilde{\alpha}}\ar[r]^-{f\Gal}&Y\ar[d]^-\beta\\
A\ar[r]^-a&B
}
$$
is commutative.
\item The correspondence sending $X\stackrel{f}{\to}Y$ to $X\Gal_f\stackrel{f\Gal}{\longrightarrow}Y$ is a functor from $\mor{G}$ to $\galmor{G}$, and this functor is left adjoint to the forgetful functor $\galmor{G}\to \mor{G}$.
\end{enumerate}
\end{mth}
\pf When $x\in X$, let $\widetilde{x}=\gamma_{X,f}(x)\in X\Gal_f$ denote its equivalence class for the relation $\sim_f$. Let $g\in G$. Then $g\widetilde{x}=\widetilde{x}$ if and only if $gx\sim_fx$, i.e. if there exists $h\in G_x^f$ such that $gx=hx$, or equivalently $h^{-1}g\in G_x$. Hence the stabilizer of $\widetilde{x}$ in $G$ is equal to $G_x^f\cdot G_x=G_x^f$, since $G_x\leq G_x^f$. So if $x,x'\in X$ are such that $f\Gal(\widetilde{x})=f\Gal(\widetilde{x}')$, i.e. if $f(x)=f(x')$, then the stabilizers $G_x^f$ of $\widetilde{x}$ and $G_{x'}^f$ of $\widetilde{x}'$ are equal, since $G_x^f$ depends only on $f(x)$. Assertion~1 follows, by Proposition~\ref{stabilizers}.\par
Let $x,x'\in X$ such that $x\sim_fx'$. Then there exist $r\in\N$, elements $z_1,\ldots, z_r$ in $X$ and elements $g_1,\ldots,g_r$ of $G$ such that $g_iz_i=z_i$ and $f(z_i)=f(x)$, for $i=1,\ldots,r$, and such that $g=g_1\cdots g_r$. It follows that
$$\beta f(z_i)=a\alpha(z_i)=\beta f(x)=a\alpha(x)\mvirg$$
for $i=1,\ldots,r$. By Proposition~\ref{stabilizers}, since $A\stackrel{a}{\to}B$ is a Galois morphism, this implies $G_{\alpha(z_i)}=G_{\alpha(x)}$. Moreover $G_z\leq G_{\alpha(z)}$ for any $z\in X$. Thus $g_i\in G_{\alpha(z_i)}=G_{\alpha(x)}$ for $i=1,\ldots,r$. It follows that $g=g_1\cdots g_r\in G_{\alpha(x)}$, hence $\alpha(x')=g\alpha(x)=\alpha(x)$. This shows the existence of a map $\widetilde{\alpha}:X\Gal_f\to A$, sending the equivalence class of $x\in X$ for $\sim_f$ to $\alpha(x)$. Such a map is obviously unique, and it is a morphism of $G$-sets. This proves Assertion~2.\par
For Assertion~3, suppose that 
$$\xymatrix{
X\ar[r]^-f\ar[d]_\alpha&Y\ar[d]^\beta\\
X'\ar[r]^{f'}&Y'
}
$$
is a morphism from $X\stackrel{f}{\to} Y$ to $X'\stackrel{f'}{\to} Y'$ in the category $\mor{G}$. Then one can compose this morphism with the morphism
$$\xymatrix{
X'\ar[r]^-{f'}\ar[d]^-{\gamma_{X',f'}}&Y'\ar@{=}[d]\\
{X'}\Gal_{f'}\ar[r]^-{{f'}\Gal}&Y'\mpoint\\
}
$$
This yields a morphism from $X\stackrel{f}{\to} Y$ to ${X'}_{f'}\Gal\stackrel{{f'}\Gal}{\longrightarrow}Y'$, which is a Galois morphism by Assertion~1. By Assertion~2, this composition factors in a unique way through the morphism $X\Gal_f\stackrel{f\Gal}{\longrightarrow}Y$. In other words there is a unique morphism of $G$-sets $\widetilde{\alpha}:X\Gal_f\to {X'}\Gal_f$ such that the diagram
$$\xymatrix{
X\ar[r]^-f\ar[drr]^(.6){\alpha}\ar[dd]_-{\gamma_{X,f}}&Y\ar[drr]^(.6){\beta}\ar@{=}[dd]&&\\
&&X'\ar[r]_{f'}\ar[dd]^-{\gamma_{X',f'}}&Y'\ar@{=}[dd]\\
{X}\Gal_f\ar[r]^-{f\Gal}\ar[drr]_(.3){\widetilde{\alpha}}&Y\ar[drr]_(.3){\beta}\\
&&{X'}\Gal_{f'}\ar[r]_-{{f'}\Gal}&Y'
}
$$
is commutative. Clearly, the map $(\alpha,\beta) \mapsto (\widetilde{\alpha},\beta)$ endows the correspondence
$$\big(X\stackrel{f}{\to}Y\big)\mapsto \big(X\Gal_f\stackrel{f\Gal}{\longrightarrow}Y\big)$$
with a structure of functor from $\mor{G}$ to $\galmor{G}$. Moreover, for any Galois morphism $A\stackrel{a}{\to}B$, Assertion~2 yields a bijection
$$\Hom_{\mor{G}}\big(X\stackrel{f}{\to}Y,A\stackrel{a}{\to}B\big)\cong \Hom_{\galmor{G}}\big(X\Gal_f\stackrel{f\Gal}{\longrightarrow}Y,A\stackrel{a}{\to}B\big)\mvirg$$
and this bijection is clearly functorial with respect to $X\stackrel{f}{\to}Y$ and $A\stackrel{a}{\to}B$. Assertion~3 follows.\findemo
\begin{rem}{Example} \label{normal closure}Let $(T,S)$ be a slice of $G$, and $f$ be the projection map from $X=G/S$ to $Y=G/T$. Then for $x=S\in X$, the group $G_x^f$ is the group generated by the stabilizers $G_{tS}$, for $t\in T$. As $G_{tS}={^tS}$, it follows that $G_x^f$ is equal to the {\em normal closure} $S^{\normal T}$ of $S$ in $T$. In this case moreover, the $G$-set $X_f\Gal$ is isomorphic to $G/(S^{\normal T})$, and the map $f\Gal$ is the projection to $G/T$.
\end{rem}

\section{The section Burnside functor}\label{The section Burnside functor}
\begin{mth}{Proposition} \label{Galois categorical product} If $X\stackrel{f}{\to}Y$ and $X'\stackrel{f'}{\to}Y'$ are Galois morphisms of $G$-sets, then 
$X\sqcup X'\stackrel{f\sqcup f'}{\longrightarrow}Y\sqcup Y'$ and $X\times X'\stackrel{f\times f'}{\longrightarrow}Y\times Y'$
are Galois morphisms of $G$-sets.
\end{mth}
\pf
The case of disjoint union follows from Proposition~\ref{base} and Corollary~\ref{restriction to subdomain}. Moreover Proposition~\ref{product} shows that if $f:X\to Y$ and $f':X'\to Y'$ are Galois morphisms of $G$-sets, then 
$$f\times f':X\times X'\to Y\times Y'$$
is a Galois morphism of $(G\times G)$-sets. By Corollary~\ref{restriction}, the restriction of this morphism to the diagonal $G\leq G\times G$ is a Galois morphism of $G$-sets. \findemo
\begin{mth}{Definition} Let $G$ be a finite group. {\em The section Burnside group $\Gamma(G)$ of $G$} is the subgroup of the slice Burnside ring $\Xi(G)$ generated by the classes of Galois morphisms of $G$-sets.\par
By Proposition~\ref{Galois categorical product}, the group $\Gamma(G)$ is actually a subring of $\Gringo(G)$, called the section Burnside ring of $G$.
\end{mth}
\begin{mth}{Lemma} \label{Galois decomposition}Let $f:X\to Y$ be a Galois morphism of finite $G$-sets. Then in the group $\Gamma(G)$
\begin{eqnarray*}
\pi(X\stackrel{f}{\to}Y)&=&\sum_{x\in [G\dom X]}\langle G_{f(x)},G_x\rangle_G\\
&=&\sum_{y\in [G\dom Y]}|G_y\dom f^{-1}(y)|\,\langle G_y,G_{x_y}\rangle_{G}\mvirg
\end{eqnarray*}
where $x_y$ is chosen in $f^{-1}(y)$, for each $y\in [G\dom f(X)]$.
\end{mth}
\pf  The first formula follows from Lemma~\ref{decomposition}. For the second one, write
$$\pi(X\stackrel{f}{\to}Y)=\sum_{y\in[G\dom Y]}\sum_{x\in [G_y\dom f^{-1}(Y)]}\pi(\gsec{G_y}{G_x}{G})\mvirg$$
and observe that $G_x=G_{x'}$ if $f(x)=f(x')$.
\findemo
\begin{mth}{Corollary} \label{Galois generators}The elements $\langle T,S\rangle_G$, where $(T,S)$ runs through a set $[\Sigma(G)]$ of representatives of conjugacy classes of sections of $G$, form a basis of $\Gamma(G)$.
\end{mth}
\pf Indeed, these elements generate $\Gamma(G)$, by Proposition~\ref{Galois decomposition}, and they are linearly independent, by~\ref{Gringo basis}.\findemo
\begin{rem}{Remark}
This also shows that $\Gamma(G)$ is the quotient of the free abelian group on the set of isomorphism classes $[X\stackrel{f}{\longrightarrow}Y]$ of {\em Galois} morphisms of finite $G$-sets, by the subgroup generated by elements of the form
$$[(X_1\sqcup X_2)\stackrel{f_1\sqcup f_2}{\longrightarrow}Y]-[X_1\stackrel{f_1}{\to}f(X_1)]-[X_2\stackrel{f_2}{\to}f(X_2)]\mvirg$$
whenever $X\stackrel{f}{\to}Y$ is a Galois morphism of finite $G$-sets with a decomposition $X=X_1\sqcup X_2$ as a disjoint union of $G$-sets, where $f_1=f_{\mid X_1}$ and $f_2=f_{\mid X_2}$ .
\end{rem}
\begin{mth}{Theorem} \begin{enumerate}
\item Let $G$ and $H$ be finite groups, and let $U$ be a finite $(H,G)$-biset. The functor
$$(X\stackrel{f}{\to}Y)\mapsto (U\times_GX\stackrel{U\times_Gf}{\longrightarrow}U\times_GY)$$
from $\galmor{G}$ to $\galmor{H}$ induces a group homomorphism 
$$\Gamma(U):\Gamma(G)\to\Gamma(H)\mpoint$$
\item The correspondence $G\mapsto \Gamma(G)$ is a Green biset functor.
\end{enumerate}
\end{mth}
\pf This follows from Theorem~\ref{Green biset functor} and Proposition~\ref{biset}. \findemo
\begin{rem}{Remark} \label{from B to Gamma} It follows from Remark~\ref{injective => galois} that the image of the morphism $i:B\to\Xi$ of Proposition~\ref{from B to Xi} is actually contained in $\Gamma$. Thus $i$ is a morphism of Green biset functors from $B$ to $\Gamma$.
\end{rem}
\section{Sections and ghost map}\label{Galois ghost}
\begin{mth}{Notation} If $(T,S)$ is a slice of $G$, denote by $\psi_{T,S}:\Gamma(G)\to \Z$ the restriction of the ring homomorphism $\phi_{T,S}:\Xi(G)\to\Z$.
\end{mth}
\begin{mth}{Lemma} \label{phi psi}Let $(T,S), (T',S')\in\Pi(G)$. Then $\psi_{T,S}=\psi_{T',S'}$ if and only if the sections $(T,S^{\normal T})$ and $(T',S'^{\normal T'})$ are conjugate in $G$ (recall that $S^{\normal T}$ denotes the normal closure of $S$ in $T$). In particular $\psi_{T,S}=\psi_{T,S^{\normal T}}$.
\end{mth}
\pf This follows from Proposition~\ref{adjunction} and Remark~\ref{normal closure}, but the following is a short direct proof~: by Lemma~\ref{phiTS}, $\psi_{T,S}=\psi_{T',S'}$ if and only if for any section $(V,U)$ of $G$
$$|\{g\in G/U\mid (T,S)^g\preceq (V,U)\}|=|\{g\in G/U\mid (T',S')^g\preceq (V,U)\}|\mpoint$$
Taking $(V,U)=(T,S^{\normal T})$ shows that there exists $g\in G$ such that $(T',S')^g\leq (T,S^{\normal T})$. This implies $(T',S'^{\normal T'})^g\preceq(T,S^{\normal T})$. Taking now $(V,U)=(T',S'^{\normal T'})$ shows that $(T',S'^{\normal T'})^{g'}\preceq(T,S^{\normal T})$, for some $g'\in G$. It follows that $(T',S'^{\normal T'})^g=(T,S^{\normal T})$.\findemo
\begin{mth}{Theorem} The map (called the {\em ghost map} for $\Gamma(G)$)
$$\Psi=\prod_{(T,S)\in[\Sigma(G)]}\psi_{T,S}:\Gamma(G)\to \prod_{(T,S)\in[\Sigma(G)]}\Z$$
is an injective ring homomorphism, with finite cokernel as morphism of abelian groups.
\end{mth}
\pf The proof is exactly the same as for Theorem~\ref{Gringo basis}~: by Lemma~\ref{generators}, the elements $\langle T,S\rangle_G$, for $(T,S)\in[\Sigma(G)]$, generate $\Gamma(G)$. Suppose that there is a non zero linear combination in the kernel of $\Psi$
$$\Lambda=\sum_{(T,S)\in[\Sigma(G)]}\lambda_{T,S}\,\langle T,S\rangle_G$$
with integer coefficients $\lambda_{T,S}\in\Z$, for $(T,S)\in[\Sigma(G)]$. Extend $\lambda$ to a function $\Sigma(G)\to \Z$, constant on conjugacy classes. Let $(Y,X)$ be an element of $\Sigma(G)$, maximal for the relation $\preceq$, such that $\lambda_{Y,X}\neq 0$. Then since by Lemma~\ref{phiTS}
$$\psi_{Y,X}(G/S\to G/T)=|\{g\in G/S\mid (Y^g,X^g)\preceq (T,S)\}|\mvirg$$
it follows that
\begin{eqnarray*}
\psi_{Y,X}(\Lambda)&=&\sum_{(T,S)\in[\Sigma(G)]}\lambda_{T,S}\psi_{Y,X}(G/S\to G/T)\\
&=&\lambda_{Y,X}\psi_{Y,X}(G/X\to G/Y)\\
&=&\lambda_{Y,X}|N_G(X,Y)/X|=0\mpoint
\end{eqnarray*}
Hence $\lambda_{Y,X}=0$, and this contradiction shows that $\Psi$ is injective (and in particular, we recover the fact that the elements $\langle T,S\rangle_G$, for $(T,S)\in[\Sigma(G)]$, form a $\Z$-basis of $\Gamma(G)$). Thus $\Psi$ is an injective morphism between free abelian groups with the same finite rank, hence it has finite cokernel.\findemo
\begin{mth}{Corollary} \label{Gamma semi-simple}Set $\Q\Gamma(G)=\Q\otimes_\Z\Gamma(G)$, and $\Q\Psi=\Q\otimes_\Z\Psi$. Then
$$\Q\Psi:\Q\Gamma(G)\to\prod_{(T,S)\in[\Sigma(G)]}\Q$$
is an isomorphism of $\Q$-algebras.
\end{mth}
\section{Sections and idempotents}\label{Sections and idempotents}
\vspace{1ex}
Corollary~\ref{Gamma semi-simple} shows that $\Q\Gamma(G)$ is a split semisimple commutative algebra. Its primitive idempotents are indexed by sections of $G$, up to conjugation~: they are the inverse images under $\Q\Psi$ of the primitive idempotents of the algebra $\mathop{\prod}_{(T,S)\in[\Sigma(G)]}\limits\Q$~:
\vspace{1ex}
\begin{mth}{Notation}
If $(T,S)\in\Sigma(G)$, denote by $\gamma_{T,S}^G$ the unique element of $\Q\otimes_\Z\Gamma(G)$ such that
$$\forall (Y,X)\in\Sigma(G),\;\;\Q\psi_{Y,X}(\gamma_{T,S}^G)=\left\{\begin{array}{cl}1&\hbox{if}\;(Y,X)=_G(T,S)\\0&\hbox{otherwise}\end{array}\right.$$
The elements $\gamma_{T,S}^G$, for $(T,S)\in\Sigma(G)$, are the primitive idempotents of $\Q\Gamma(G)$.
\end{mth}
\vspace{4ex}
\pagebreak[4]
\begin{mth}{Theorem} \label{galois idempotent}Let $\Lln$ denote the restriction of the relation $\preceq$ to $\Sigma(G)$. Then for $(T,S)\in \Sigma(G)$
$$\gamma_{T,S}^G=\frac{1}{|N_G(T,S)|}\sum_{(V,U)\lln (T,S)}|U|\mu_\Sigma\big((V,U),(T,S)\big)\gsect{V,U}_G\mvirg$$
where $\mu_\Sigma$ is the M\"obius function of the poset $(\Sigma(G),\Lln)$.
\end{mth}
\pf The proof is the same as the proof of Theorem~\ref{slices and idempotents}.\findemo
\begin{mth}{Proposition} \label{mobius varphi}Let $(\mathcal{X},\leq)$ be a finite poset. Let $\varphi:\mathcal{X}\to\mathcal{X}$ be a map of posets such that $\varphi\circ\varphi=\varphi$ and $x\leq \varphi(x)$ for any $x\in\mathcal{X}$. Then the M\"obius function $\mu_\mathcal{Y}$ of the subposet $\mathcal{Y}=\varphi(\mathcal{X})$ of $\mathcal{X}$ is given by
\begin{equation}\label{mobius phi}
\forall y,z\in\mathcal{Y},\;\;\mu_\mathcal{Y}(y,z)=\sumb{y\leq u\in\mathcal{X}}{\rule{0ex}{1.2ex}\varphi(u)=z}\mu_{\mathcal{X}}(y,u)\mvirg
\end{equation}
where $\mu_{\mathcal{X}}$ is the M\"obius function of $(\mathcal{X},\leq)$.
\end{mth}
\pf For $y,z\in\mathcal{Y}$, denote by $m(y,z)$ the right hand side of Equation~\ref{mobius phi}. Then for $y,t\in\mathcal{Y}$
\begin{eqnarray*}
\sumb{z\in\mathcal{Y}}{y\leq z\leq t}m(y,z)&=&\sumb{z\in\mathcal{Y}}{y\leq z\leq t}\sumb{y\leq u\in\mathcal{X}}{\varphi(u)=z}\mu_\mathcal{X}(y,u)\\
&=&\sum_{(z,u)\in\mathcal{P}}\mu_\mathcal{X}(y,u)\mvirg
\end{eqnarray*}
where 
$$\mathcal{P}=\{(z,u)\mid z\in\mathcal{Y},\,u\in\mathcal{X},\,y\leq u\leq\varphi(u)=z\leq t\}\mpoint$$
Set $\mathcal{Q}=\{u\in\mathcal{X}\mid y\leq u\leq t\}$. If $(z,u)\in\mathcal{P}$, then clearly $u\in\mathcal{Q}$. Conversely, if $u\in\mathcal{Q}$, then $\big(\varphi(u),u\big)\in\mathcal{P}$~: indeed $\varphi(u)\in \mathcal{Y}=\varphi(\mathcal{X})$, and moreover $\varphi(u)\leq\varphi(t)=t$, since $t\in\varphi(\mathcal{X})$ and $\varphi\circ\varphi=\varphi$. Now the maps $(z,u)\in\mathcal{P}\mapsto u\in\mathcal{Q}$ and $u\in\mathcal{Q}\mapsto \big(\varphi(u),u\big)\in\mathcal{P}$ are mutual inverse bijections. It follows that
$$\sumb{z\in\mathcal{Y}}{\rule{0ex}{1.2ex}y\leq z\leq t}m(y,z)=\sumb{u\in\mathcal{X}}{\rule{0ex}{1.2ex}y\leq u\leq t}\mu_\mathcal{X}(y,u)\mpoint$$
This is equal to $1$ if $y=t$, and to zero otherwise. The proposition follows.\findemo
\begin{mth}{Corollary} \label{Galois idempotents}Let $(V,U)$ and $(T,S)$ be sections of $G$. Then
\begin{equation}\label{mobius sections}\mu_{\Sigma(G)}\big((V,U),(T,S)\big)=\mu(V,T)\Big(\sumb{U\leq X\leq V}{\rule{0ex}{1.7ex}X^{\normal T}=S}\mu(U,X)\Big)\mpoint
\end{equation}
In particular in $\Q\Gamma(G)$
$$\gamma_{T,S}^G=\frac{1}{|N_G(T,S)|}\sumc{U\normal V\leq T}{\rule{0ex}{1.4ex}U\leq X\leq V}{\rule{0ex}{1.7ex}X^{\normal T}=S}|U|\mu(U,X)\mu(V,T)\,\gsect{V,U}_G\mpoint$$
\end{mth}
\pf For Equation~\ref{mobius sections}, apply Proposition~\ref{mobius varphi} to the poset $\mathcal{X}=\big(\Pi(G),\preceq\nolinebreak\big)$, and to the map $\varphi:\mathcal{X}\to \mathcal{X}$ defined by
$$\varphi\big((Y,X)\big)=(Y,X^{\normal Y})\mvirg$$
and then use Equation~\ref{poset Pi}. Then substitute the value of $\mu_{\Sigma}\big((V,U),(T,S)\big)$ in the formula of Theorem~\ref{galois idempotent}.\findemo 
\section{The image of the ghost map}\label{Galois ghost map image}
The following is a characterization of the image of the ghost map for the section Burnside ring, similar to Theorem~\ref{ghost image}~:
\begin{mth}{Theorem} \label{Galois ghost image} Let $G$ be a finite group. Let $\mathsf{m}=\big(m({T,S})\big)_{(T,S)\in\Sigma(G)}$ be a sequence of integers indexed by $\Sigma(G)$, constant on $G$-conjugacy classes of sections. Then the sequence $[\mathsf{m}]=\big(m({T,S})\big)_{(T,S)\in[\Sigma(G)]}$ of representatives lies in the image of the ghost map~$\Psi$ if and only if, for any section $(T,S)$ of~$G$
$$\sum_{g\in N_G(T,S)/S}m\big({\gen{gT},\gen{gS}^{\normal{\gen{gT}}}}\big)\equiv 0\;\big({\rm mod.}\,|N_G(T,S)/S|\big)\mpoint$$
\end{mth}
\pf The proof is very similar to the proof of Theorem~\ref{ghost image}, with two differences~: the first one is that the poset $\Pi(G)$ of slices of $G$ has to be replaced by the poset $\Sigma(G)$ of sections of $G$. The second one is that the slice $(\gen{gT},\gen{gS})$ has to be changed to the corresponding section $\big({\gen{gT},\gen{gS}^{\normal{gT}}}\big)$ of $G$. Apart from these two differences, the proof goes through without changes.\findemo
\section{Prime spectrum}\label{Galois prime spectrum}
\begin{mth}{Notation} Let $p$ denote either 0 or a prime number.
\begin{itemize}
\item If $(T,S)\in\Sigma(G)$, denote by  $J_{T,S,p}$ the prime ideal $I_{T,S,p}\cap \Gamma(G)$ of $\Gamma(G)$. In other words $J_{T,S,p}$ is the kernel of the ring homomorphism 
$$\Gamma(G)\stackrel{\psi_{T,S}}{\longrightarrow}\Z\to \Z/p\Z\mvirg$$
where the right hand side map is the projection.
\item Let $\Omega(G)$ denote the set of triples $(T,S,p)$, where $(T,S)\in\Sigma(G)$ is such that $|N_G(T,S)/S|\not\equiv 0 \;({\rm mod.}\;p)$.
\end{itemize}
\end{mth}
The group $G$ acts on $\Omega(G)$, by $^g(T,S,p)=({^gT},{^gS},p)$, for $g\in G$, and the ideal $J_{T,S,p}$ only depends on the $G$-orbit of $(T,S,p)$. Conversely~:

\begin{mth}{Proposition} Let $I$ be a prime ideal of $\Gamma(G)$, and $R=\Gamma(G)/I$. Denote by $\psi:\Gamma(G)\to R$ the projection map, and denote by $p\geq 0$ the characteristic of $R$. Then $R\cong \Z/p\Z$ and~:
\begin{enumerate}
\item If $p=0$, there exists a section $(T,S)$ of $G$ such that $\psi=\psi_{T,S}$, and $(T,S)$ is unique up to $G$-conjugation, with this property.
\item If $p>0$, there exists a section $(T,S)$ of $G$ such that $\psi$ is the reduction modulo $p$ of $\psi_{T,S}$ and $N_G(T,S)/S$ is a $p'$-group, and $(T,S)$ is unique up to $G$-conjugation, with these properties.
\end{enumerate}
In particular, there exists a unique $(T,S,p)\in\Omega(G)$, up to conjugation, such that $I=J_{T,S,p}$.
\end{mth}
\pf The proof is exactly the same as the proof of Proposition~\ref{prime ideals}, with slices replaced by sections~: let $(T,S)$ be a section of $G$, minimal for the relation $\preceq$, such that $\gsect{T,S}_G\notin I$. Then by Proposition~\ref{product formula}, for any $(Y,X)\in\Sigma(G)$
\begin{eqnarray*}
\gsect{T,S}_G\gsect{Y,X}_G&=&\sum_{g\in[S\dom G/X]}\gsect{T\cap {^gY},S\cap{^gX}}_G\\
&\equiv&\sumc{g\in G/X}{S\leq{^gX}}{T\leq{^gY}}\gsect{T,S}_G \;\;({\rm mod.}\, I)\\
&=&\psi_{T,S}\big(\gsect{Y,X}_G\big)\gsect{T,S}_G\mpoint
\end{eqnarray*}
Since $I$ is prime, it follows that $\gsect{Y,X}_G-\psi_{T,S}(\gsect{Y,X}_G)1_{\Gamma(G)}\in I$. In particular $R=\Gamma(G)/I$ is generated by the image of $1_{\Gamma(G)}$, hence $R\cong \Z/p\Z$, where $p$ is the characteristic of $R$. Since $R$ is an integral domain, the number $p$ is either 0 or a prime.\begin{enumerate}
\item If $p=0$, then $R=\Z$, and $\psi=\psi_{T,S}$. And if $(T',S')\in\Sigma(G)$ is such that $\psi_{T,S}=\psi_{T',S'}$, then both $\psi_{T,S}\big(\langle T',S'\rangle_G\big)$ and $\psi_{T',S'}\big(\langle T,S\rangle_G\big)$ are non zero. Then there exist elements $g,g'\in G$ such that $(T^g,S^g)\preceq (T',S')$ and $(T'^{g'},S'^{g'})\preceq (T,S)$, so $(T,S)$ and $(T',S')$ are conjugate in $G$.
\item If $p>0$, then $R=\Z/p\Z$, and $\psi$ is equal to the reduction of $\psi_{T,S}$ modulo~$p$. Since $\psi\big(\gsect{T,S}_G\big)=|N_G(T,S)/S|$ is non zero in $R$, it follows that $N_G(T,S)/S$ is a $p'$-group. If $(T',S')$ is another section of $G$ such that $\psi$ is the reduction modulo $p$ of $\psi_{T',S'}$, and $N_G(T',S')/S'$ is a $p'$-group, then 
$$|N_G(T,S)/S|=\psi_{T,S}\big(\gsect{T,S}_G\big)\equiv \psi_{T',S'}\big(\gsect{T,S}_G\big)\;\;({\rm mod.}\,p)\mpoint$$ 
This is non zero. Similarly $|N_G(T,S)/S|\equiv \psi_{T,S}\big(\gsect{T',S'}_G\big)\;\;({\rm mod.}\,p)$ is non zero. In particular $\psi_{T',S'}\big(\gsect{T,S}_G\big)$ and $\psi_{T,S}\big(\gsect{T',S'}_G\big)$ are both non zero, and it follows as above that $(T,S)$ and $(T',S')$ are conjugate in~$G$.\findemo
\end{enumerate}
\begin{mth}{Notation} Let $p$ be a prime number. 
\begin{itemize}
\item Let $\Sigma_p(G)$ denote the set of sections $(T,S)$ of $G$ such that $N_G(T,S)/S$ is a $p'$-group. In other words $\Sigma_p(G)=\Sigma(G)\cap\Pi_p(G)$.
\item For any $(T,S)\in\Sigma(G)$, let $(T,S)\compn{p}$ denote the unique element $(V,U)$ of $\Sigma_p(G)$, up to conjugation, such that $J_{T,S,p}=J_{V,U,p}$.
\end{itemize}
\end{mth}
\begin{mth}{Proposition} Let $p$ be a prime number. If $(T,S)$ is a section of $G$, let $(T,S)\plusn{p}$ denote a section of the form $\big(PT,(PS)^{\normal PT}\big)$ of $G$, where $P$ is a Sylow $p$-subgroup of $N_G(T,S)$. \par
Define inductively an increasing sequence $(T_n,S_n)$ in $(\Sigma(G),\preceq)$ by $(T_0,S_0)=(T,S)$, and $(T_{n+1},S_{n+1})=(T_n,S_n)\plusn{p}$, for $n\in\N$. Then $(T,S)\compn{p}$ is conjugate to the largest term $(T_\infty,S_\infty)$ of the sequence $(T_n,S_n)$. 
\end{mth}
\pf Again, the proof is the same as the proof of Proposition~\ref{closure}~: the section $(T,S)\compn{p}$ is a minimal element $(V,U)$ of the poset $(\Sigma(G),\preceq)$ such that
$$\psi_{T,S}(V,U)=|\{g\in G/U\mid (T^g,S^g)\preceq (V,U)\}|\not\equiv 0\;({\rm mod.} p)\mpoint$$
Thus one can assume that $(T,S)\preceq (V,U)$.  But $\psi_{T,S}\equiv \psi_{PT,(PS)^{\normal{PT}}}\;({\rm mod.}\;p)$ by Corollary~\ref{modp} and Lemma~\ref{phi psi}, for any $p$-subgroup $P$ of $N_G(T,S)$, hence one can also assume that $(T,S)\plusn{p}\preceq (V,U)$, hence that $(T_\infty,S_\infty)\preceq (V,U)$, by induction. Moreover $\psi_{T_\infty,S_\infty}\equiv \psi_{V,U}\;({\rm mod.}\;p)$. As $N_G(T_\infty,S_\infty)/S_\infty$ is a $p'$-group, it follows that $(T_\infty,S_\infty)= (V,U)$, as was to be shown. \findemo
\pagebreak[3]
\begin{mth}{Proposition} \label{galois ideal inclusion}Let $(T,S,p)$, $(T',S',p')$ be elements of $\Omega(G)$. Then $J_{T',S',p'}\subseteq J_{T,S,p}$ if and only if
\begin{itemize}
\item either $p'=p$ and the sections $(T',S')$ and $(T,S)$ are conjugate in~$G$.
\item or $p'=0$ and $p>0$, and the sections $(T',S')\compn{p}$ and $(T,S)$ are conjugate in~$G$.
\end{itemize}
\end{mth}
\pf (see Proposition~\ref{ideal inclusion}) Set $J=J_{T,S,p}$ and $J'=J_{T',S',p'}$. Then $\Gamma(G)/J'\cong \Z/p'\Z$ maps surjectively to $\Gamma(G)/J\cong\Z/p\Z$. Thus if $p=p'$, this projection map is an isomorphism, hence $J=J'$ and the sections $(T,S)$ and $(T',S')$ are conjugate in~$G$. And if $p\neq p'$, then $p'=0$ and $p>0$. The morphism $\psi_{T,S}$ is equal to the reduction modulo $p$ of the morphism $\psi_{T',S'}$. In other words $J_{T',S',p}=J_{T,S,p}$, hence $(T,S)$ is conjugate to $(T',S')\compn{p}$.\findemo
\begin{mth}{Corollary} Let $p$ be a prime number, and let $\Z_{(p)}$ be the localization of $\Z$ at the set $\Z-p\Z$. Let $\Omega_p(G)$ denote the subset of $\Omega(G)$ consisting of triples $(T,S,0)$, for $(T,S)\in\Sigma(G)$, and $(T,S,p)$, for $(T,S)\in\Sigma_p(G)$. Then~:
\begin{enumerate}
\item The prime ideals of the ring $\Z_{(p)}\Gamma(G)$ are the ideals $\Z_{(p)}J_{T,S,q}$, for $(T,S,q)\in\Omega_p(G)$.
\item If $(T,S,q), (T',S',q')\in\Theta_p(G)$, then $\Z_{(p)}J_{T',S',q'}\subseteq \Z_{(p)}J_{T,S,q}$ if and only~if~:
\begin{itemize}
\item either $q=q'$, and the sections $(T,S)$ and $(T',S')$ are conjugate in~$G$.
\item or $q'=0$, $q=p$, and the sections $(T',S')\compn{p}$ and $(T,S)$ are conjugate in $G$.
\end{itemize}
\item The connected components of the spectrum of $\Z_{(p)}\Gamma(G)$ are indexed by the conjugacy classes of $\Sigma_p(G)$. The component indexed by $(T,S)\in\Sigma_p(G)$ consists of a unique maximal element $\Z_{(p)}J_{T,S,p}$, and of the ideals $\Z_{(p)}J_{T',S',0}$, where $(T',S')\in\Sigma(G)$ is such that $(T',S')\compn{p}$ is conjugate to $(T,S)$ in $G$.
\end{enumerate}
\end{mth} 
\pf (See Corollary~\ref{Zp components}) The prime ideals of $\Z_{(p)}\Gamma(G)$ are of the form $\Z_{(p)}I$, where $I$ is a prime ideal of $\Gamma(G)$ such that $I\cap(\Z-p\Z)=\emptyset$. Equivalently $I=J_{T,S,0}$ or $I=J_{T,S,p}$. This proves Assertion~1. Now Assertion~2 follows from Proposition~\ref{galois ideal inclusion}, and Assertion~3 follows from Assertion~2.\findemo
\begin{mth}{Corollary} \label{galois pi idempotents}\begin{enumerate}
\item The primitive idempotents of the ring $\Z_{(p)}\Gamma(G)$ are indexed by the conjugacy classes of $\Sigma_p(G)$. The primitive idempotent $\varepsilon_{V,U}^G$ indexed by $(V,U)\in\Sigma_p(G)$ is equal to
$$\varepsilon_{V,U}^G=\sumb{(T,S)\in[\Sigma(G)]}{(T,S)\compni{p}=_G(T,S)}\gamma_{T,S}^G\mpoint$$
\item Let $\pi$ be a set of prime numbers, and $\Z_{(\pi)}$ be the localization of $\Z$ relative to $\Z-\cup_{p\in\pi}p\Z$. Let $\mathcal{F}$ be a set of sections of $G$, invariant by $G$-conjugation, and $[\mathcal{F}]$ be a set of representatives of $G$-conjugacy classes of $\mathcal{F}$. Then the following conditions are equivalent~:
\begin{enumerate}
\item The idempotent 
$$\gamma_\mathcal{F}^G=\sum_{(T,S)\in[\mathcal{F}]}\gamma_{T,S}^G$$
of $\Q\Gamma(G)$ lies in $\Z_{(\pi)}\Gamma(G)$.
\item Let $(T,S)\in \Sigma(G)$, and let $P$ be a $p$-subgroup of $N_G(T,S)$, for some $p\in\pi$. Then $(T,S)\in\mathcal{F}$ if and only if $\big(PT,(PS)^{\normal PT}\big)\in\mathcal{F}$. 
\end{enumerate}
\end{enumerate}
\end{mth}
\pf The proof is almost identical to the proof of Corollary~\ref{pi idempotents}~: let $\mathcal{F}$ be a set of sections of $G$, invariant by $G$-conjugation, and $[\mathcal{F}]$ be a set of representatives of $G$-conjugacy classes of $\mathcal{F}$. The idempotent
$$\gamma_\mathcal{F}^G=\sum_{(T,S)\in[\mathcal{F}]}\gamma_{T,S}^G$$
of $\Q\Gamma(G)$ lies in $\Z_{(p)}\Gamma(G)$, for some prime $p$, if and only if there exists an integer $m$, not divisible by $p$, such that $u=m\gamma_\mathcal{F}^G\in\Gamma(G)$. Let $(T,S)\in\Sigma(G)$, and let $P$ be a $p$-subgroup of $N_G(T,S)$. The integer $\psi_{T,S}(u)$ is equal to $m$ if $(T,S)\in\mathcal{F}$, and to 0 otherwise. Hence it is coprime to $p$ if and only if $(T,S)\in\mathcal{F}$. Since $\psi_{T,S}$ and $\psi_{PT,(PS)^{\normal PT}}$ are congruent modulo $p$, it follows that $(T,S)\in\mathcal{F}$ if and only if $\big(PT,(PS)^{\normal PT}\big)\in\mathcal{F}$. \par
Hence if $(T,S)$ and $(T',S')$ are sections of $G$ such that $(T,S)\compn{p}=_G(T',S')\compn{p}$, then $(T,S)\in\mathcal{F}$ if and only if $(T',S')\in\mathcal{F}$. Thus $\mathcal{F}$ is a disjoint union of sets of the form 
$$E_{V,U}=\{(T,S)\in\Sigma(G)\mid (T,S)\compn{p}=_G(V,U)\}\mvirg$$
for some sections $(V,U)\in\Sigma_p(G)$. In other words the idempotent $\gamma_\mathcal{F}^G$ is a sum of some idempotents $\varepsilon_{V,U}^G$, for $(V,U)\in\Sigma_p(G)$. \par
But the primitive idempotents of the ring $\Z_{(p)}\Gamma(G)$ are in one to one correspondence with the connected components of its spectrum, which precisely are indexed by the conjugacy classes of $\Sigma_p(G)$. It follows that $\gamma_\mathcal{F}^G=\varepsilon_{V,U}^G$ is equal to the idempotent corresponding to the component indexed by $(V,U)$, for any $(V,U)\in\Sigma_p(G)$. This proves Assertion~1. This also proves Assertion~2 in the case where $\pi$ consists of a single prime number.\par
For the general case, observe that $\gamma_{\mathcal{F}}^G$ lies in $\Z_{(\pi)}\Gamma(G)$ if and only if it lies in $\Z_{(p)}\Gamma(G)$, for any $p\in \pi$.\findemo
The following proposition is the analogue of Theorem~\ref{connected Xi}, for the ring $\Gamma(G)$. Its statement is a bit simpler~:
\begin{mth}{Theorem} \label{Galois connected}Let $G$ be a finite group. 
\begin{enumerate}
\item The primitive idempotents of $\Gamma(G)$ are indexed by the conjugacy classes of perfect subgroups of $G$. The idempotent $\gamma_H^G$ indexed by the perfect subgroup $H$ is equal to
$$\gamma_H^G=\sumb{(T,S)\in[\Sigma(G)]}{D^\infty(T)=_GH}\gamma_{T,S}^G\mvirg$$
where $D^\infty(T)$ denotes the last term in the derived series of $T$.
\item The prime spectrum of $\Gamma(G)$ is connected if and only if $G$ is solvable.
\end{enumerate}
\end{mth}
\pf As in the proof of Theorem~\ref{connected Xi}, let $\sim$ denote the finest equivalence relation on the set $\Sigma(G)$ such that for any $(T,S), (T',S')\in\Sigma(G)$
$$\exists p,\;(T,S)\compn{p}=_G(T',S')\compn{p}\implies (T,S)\sim(T',S')\mpoint$$
If we can show that 
$$(T,S)\sim (T',S')\Leftrightarrow D^\infty(T)=_GD^\infty(T')\mvirg$$
then Assertion 1 follows from Corollary~\ref{galois pi idempotents}, applied to the set $\pi$ of all primes.\par
Clearly, if there exists a $p$-subgroup $P\leq N_G(T,S)$ such that $(T',S')$ is conjugate to $\big(PT,(PS)^{\normal PT}\big)$, then $T'$ is conjugate to $PT$, and $D^\infty(T')$ is conjugate to $D^\infty(T)$. By transitivity, for any $(T,S), (T',S')\in\Sigma(G)$, if $(T,S)\sim (T',S')$, then $D^\infty(T)=_GD^\infty(T')$.\par
To show the converse, is it enough to show that $(T,S)\sim\big(D^\infty(T),D^\infty(T)\big)$, for any $(T,S)\in\Sigma(G)$. Let $p$ be any prime, and $P$ be a Sylow $p$-subgroup of $T$. Then $(T,S)\compn{p}=\big(T,(PS)^{\normal PT}\big)\compn{p}$, hence $(T,S)\sim (T,S')$, where $S'=(PS)^{\normal PT}$ is a normal subgroup of $T$, containing $S$, and of $p'$-index in $T$. Since $p$ was arbitrary, it follows by induction that $(T,S)\sim (T,T)$. 
Now for any prime $p$, if $P$ is a Sylow $p$-subgroup of $T$, then $P\leq N_G\big(O^p(T)\big)$, and $T=PO^p(T)$. It follows that $(T,T)\compn{p}=\big(O^p(T),O^p(T)\big)\compn{p}$. Again, since $p$ is arbitrary, it follows that $(T,T)\sim\big(D^\infty(T),D^\infty(T)\big)$, and this completes the proof of Assertion~1.\par
Assertion~2 follows easily, since by Assertion~1, the spectrum of $\Gamma(G)$ is connected if and only if the trivial group is the only perfect subgroup of $G$, i.e. if $G$ is solvable.\findemo
\begin{mth}{Corollary} The images of the primitive idempotents of $B(G)$ by the morphism $i_G:B(G)\to\Gamma(G)$ are the primitive idempotents of $\Gamma(G)$. In other words, if $H$ is a perfect subgroup of $G$, then
$$\gamma_H^G=\frac{1}{|N_G(H)|}\sum_{K\leq H}|K|\mu(K,H)\gsect{K,K}_G\mpoint$$
\end{mth}
\pf Indeed by a theorem of Dress (\cite{dressresoluble}), the primitive idempotents of $B(G)$ are indexed by the conjugacy classes of perfect subgroups of $G$. As the morphism $i_G:B(G)\to \Gamma(G)$ is an injective unital ring homomorphism (see Remark~\ref{from B to Gamma} and Proposition~\ref{from B to Xi}), it follows that the primitive idempotents of $B(G)$ are mapped to primitive idempotents in $\Gamma(G)$, and that every primitive idempotent of $\Gamma(G)$ is in the image of $B(G)$.\findemo
\section{Unit group}\label{Galois units}
All the results of Section~\ref{units} about the unit group of the slice Burnside ring have an analogue for the group $\Gamma(G)^\times$ of the section Burnside ring of a finite group $G$. Namely~:\spn
$\bullet$ The restricted ghost map yields an injective group homomorphism 
$$\Psi^\times:\Gamma(G)^\times\hookrightarrow \prod_{(T,S)\in\Sigma(G)}\Z^\times\mpoint$$
The following lemma follows~:
\begin{mth}{Lemma} Let $G$ be a finite group, and let $u\in\Gamma(G)$. The following conditions are equivalent~:
\begin{enumerate}
\item $u\in \Gamma(G)^\times$.
\item $\psi_{T,S}(u)\in\{\pm 1\}$, for any $(T,S)\in\Sigma(G)$.
\item $u^2=1$.
\end{enumerate}
In particular $\Gamma(G)^\times$ is a finite elementary abelian 2-group.\par
\end{mth}
\spn$\bullet$ The following is an analogue of Proposition~\ref{Feit-Thompson}~:
\begin{mth}{Proposition} Feit-Thompson's theorem is equivalent to the statement that, if $G$ has odd order, then $\Gamma(G)^\times=\{\pm 1\}$.
\end{mth}
\pf The proof is the same as the proof of Proposition~\ref{Feit-Thompson}~: the argument uses only the formulae for the primitive idempotents of $\Q \Gamma(G)$ (Corollary~\ref{Galois idempotents}), and the characterization of solvable groups by the connectedness of the prime spectrum of $\Gamma(G)$ (Proposition ~\ref{Galois connected}).\findemo
\spn$\bullet$ The following theorem is an analogue of Yoshida's characterization (\cite{yoshidaunit}) of the unit group of the usual Burnside ring~:
\begin{mth}{Theorem} \label{Galois units characterization}Let $G$ be a finite group, and let $\mathsf{m}=\big(m(T,S)\big)_{(T,S)\in\Sigma(G)}$ be a sequence of integers in $\{\pm 1\}$ indexed by $\Sigma(G)$, constant on $G$-conjugacy classes of sections. Then the sequence $[m]= \big(m(T,S)\big)_{(T,S)\in[\Sigma(G)]}$ of representatives lies in the image of the restricted ghost map $\Psi^\times$ if and only if for any $(T,S)\in\Sigma(G)$, the map
$$g\in N_G(T,S)/S\mapsto m\big(\gen{gT},\gen{gS}^{\normal\gen{gT}}\big)/m(T,S)\in\{\pm 1\}$$
is a group homomorphism.
\end{mth}
\pf Again, the proof is the same as for Theorem~\ref{units characterization}~: it only requires Theorem~\ref{Galois ghost image} and Lemma~\ref{phi psi}.\findemo
\spn $\bullet$ Finally, the correspondence $G\mapsto \Gamma(G)^\times$ is not a biset functor~:
\begin{mth}{Proposition} \label{Galois not biset functor}The correspondence sending a finite group $G$ to $\Gamma(G)^\times$ cannot be endowed with a structure of biset functor.
\end{mth}
\pf The proof of Proposition~\ref{not biset functor} goes through without change here~: the argument uses computations for the trivial group, the group $C_2$, and the group $(C_2)^2$. All these groups are abelian, hence sections and slices coincide.\findemo
\vspace{4ex}\par
\pagebreak[3]
\centerline{\Large\bf Appendix}\vspace{4ex}\par\refstepcounter{section}
\npar Let $G$ and $H$ be finite groups, and $U$ be a finite $(H,G)$-biset. Let $U\op$ denote the {\em opposite biset}, i.e. the $(G,H)$-biset equal to $U$ as a set, with actions reversed by taking inverses (i.e. $g\cdot u\cdot h\;\hbox{[in $U\op$]}=h^{-1}ug^{-1}\;\hbox{[in $U$]}$, for $g\in G$, $u\in U$, and $h\in H$). \par
When $X$ is a finite $G$-set, the set $\Hom_{G\hbox{-}\mathsf{set}}(U\op,X)$ is a finite $H$-set~: if $\varphi:U\op\to X$ is a morphism of $G$-sets, and if $h\in H$, then the morphism $h\varphi:U\op\to X$ is defined by $(h\varphi)(u)=\varphi(h^{-1}u)$, for $u\in U$. This correspondence
$$T_U:X\mapsto \Hom_{G\hbox{-}\mathsf{set}}(U\op,X)$$
is actually a functor from the category $\gset{G}$ of finite $G$-sets to the category $\gset{H}$. One can show (see e.g. Section 11.2.13 of \cite{bisetfunctors}) that this functor induces a map $t_U:B(G)\to B(H)$ between the usual Burnside rings of $G$ and $H$, called {\em the generalized tensor induction} with respect to $U$. This induction is not additive, but multiplicative (i.e. $t_U(ab)=t_U(a)t_U(b)$, for any $a,b\in B(G)$). It yields a biset functor structure on the unit group of the usual Burnside ring.\par
A natural question is to know whether this construction can be extended to the rings $\Xi(G)$ and $\Gamma(G)$~: indeed, if $X\stackrel{f}{\to}Y$ is a morphism of finite $G$-sets, then the morphism $T_U(X\stackrel{f}{\to}Y)$
$$\xymatrix{
\Hom_{G\hbox{-}\mathsf{set}}(U\op,X)\ar[rrr]^-{\Hom_{G\hbox{-}\mathsf{set}}(U\op,f)}&&&\Hom_{G\hbox{-}\mathsf{set}}(U\op,Y)
}
$$
is a morphism of finite $H$-sets. Does this correspondence induce a map $\Xi(G)\to \Xi(H)$ or a map $\Gamma(G)\to\Gamma(H)$~?\par
In other words~:
\begin{itemize}
\item[(Q1)] Are the defining relations of $\Xi(G)$ preserved by $T_U$~?
\item[(Q2)] Does $T_U$ map a Galois morphism to a Galois morphism~? \par
\end{itemize}
The answer to these two questions has to be no in general, for otherwise $T_U$ would yield a biset functor structure on the unit groups of the slice Burnside ring and the section Burnside ring, contradicting Proposition~\ref{not biset functor} and Proposition~\ref{Galois not biset functor}. But the answer is yes with the additional assumption that the biset $U$ be {\em left inert}, according to the following definition~:
\begin{mth}{Definition} An $(H,G)$-biset $U$ is called {\em left inert} if $Hu\subseteq uG$, for any $u\in U$, in other words if $H$ acts trivially on the set of orbits $U/G$.
\end{mth}
\begin{rem}{Example} \label{example left inert}\begin{itemize}
\item If $U$ is transitive as a right $G$-set, then $U$ is left inert, since $|U/G|=1$. In particular, the identity $(G,G)$-biset $G$ is left inert.  Conversely, if $U$ is left inert, then each right orbit $uG$, for $u\in U$, is an $(H,G)$-biset, so $U$ is a disjoint union of $(H,G)$-bisets which are transitive as right $G$-sets.
\item The disjoint union of two left inert $(H,G)$-bisets is left inert. Any sub-biset of a left inert biset if left inert.
\item Left inert bisets can be composed~: if $G$, $H$, and $K$ are groups, if $U$ is a left inert $(H,G)$-biset, and if $V$ is a left inert $(K,H)$-biset, then $V\times_HU$ is a left inert $(K,G)$-biset~: indeed, for $u\in U$ and $v\in V$
$$K(v,_{_H}u)=(Kv,_{_H}u)\subseteq (vH,_{_H}u)=(v,_{_H}Hu)\subseteq (v,_{_H}uG)=(v,_{_H}u)G\mvirg$$
where $(v,_{_H}u)$ denotes the image of $(v,u)$ in $V\times_HU$.
\end{itemize}
\end{rem}
\npar The following proposition deals with Question (Q2) above, in the case of a left inert biset~:
\begin{mth}{Proposition} \label{Galois left inert}Let $X\stackrel{f}{\to} Y$ be a Galois morphism of $G$-sets. If~$U$ is a left inert $(H,G)$-biset, then the morphism $T_U(X\stackrel{f}{\to} Y)$ is a Galois morphism of $H$-sets.
\end{mth}
\pf Let $\varphi,\psi\in T_U(X)=\Hom_{G\hbox{-}\mathsf{set}}(U\op,X)$ having the same image in $T_U(Y)$. This means that $f\circ \varphi=f\circ\psi$, i.e. that $f\big(\varphi(u)\big)=f\big(\psi(u)\big)$, for any $u\in U$. Since $f$ is a Galois morphism, it follows that $G_{\varphi(u)}=G_{\psi(u)}$, for any $u\in U$. \par
Suppose that $h\in H$ stabilizes $\varphi$. It means that $\varphi(h^{-1}u)=\varphi(u)$, for any $u\in U$. But since $Hu\subseteq uG$, there exists $g\in G$ (depending on $u$) such that $hu=ug$, i.e. $h^{-1}u=ug^{-1}$. Then 
$$\varphi(u)=\varphi(h^{-1}u)=\varphi(ug^{-1})=g\varphi(u)\mvirg$$
hence $g\in G_{\varphi(u)}$. It follows that
$$\psi(u)=g\psi(u)=\psi(ug^{-1})=\psi(h^{-1}u)=(h\psi)(u)\mpoint$$
Hence $h$ stabilizes $\psi$, and by symmetry $\varphi$ and $\psi$ have the same stabilizer in~$H$. This shows that $T_U(X\stackrel{f}{\to} Y)$ is a Galois morphism of $H$-sets.\findemo
\begin{rem}{Remark} \label{counter example Galois}The following example shows that Lemma~\ref{Galois left inert} is no longer true without the hypothesis that $U$ is left inert~: suppose that $G=\un$, and that $H$ is non trivial. Let $U=H$, viewed as an $(H,G)$-biset. Then $U$ is not left inert. Let $X$ be a set of cardinality 2, let $Y$ be a set of cardinality~1, and let $f:X\to Y$ be the unique map. Then $f$ is a Galois morphism of $G$-sets, for trivial reasons.\par
In this case $\Hom_{G}(U\op,X)$ is isomorphic to the set ${\bf 2}^H$ of subsets of $H$, on which $H$ acts by left translation. The set $\Hom_{G}(U\op,Y)$ is a set $\bullet$ of cardinality 1, and the map $T_U(f)$ is the only possible map ${\bf 2}^H\to\bullet$. In particular, all the elements of ${\bf 2}^H$ have the same image. But the subset $\{1\}$ of $H$ has a trivial stabilizer, whereas the subset $H$ of $H$ has stabilizer equal to~$H$. Thus $T_U(X\stackrel{f}{\to}Y)$ is not a Galois morphism of $H$-sets, by Proposition~\ref{stabilizers}.
\end{rem}
\begin{rem}{Remark} The example given in Remark~\ref{counter example Galois} also shows that the answer to Question (Q1) is no, in general~: indeed, keeping the same notation, the image of the morphism $X\stackrel{f}{\to}Y$ in $\Xi(G)$ is equal to the sum of two copies of the image of $Y\stackrel{\Id}{\to} Y$. In other words, in $\Xi(G)\cong\Z$,
$$\pi(X\stackrel{f}{\to}Y)=2=\pi(X\stackrel{\Id}{\to} X)\mpoint$$
But $T_U(X\stackrel{f}{\to}Y)$ is isomorphic to ${\bf 2}^H\to\bullet$, hence
\begin{equation}\label{sum example}\pi\big(T_U(X\stackrel{f}{\to}Y)\big)=\sumb{A\subseteq H}{A\,{\rm mod. H}}\gsect{H,H_A}_H\end{equation}
by Lemma~\ref{decomposition}, where the summation runs over a set of representatives of {\em subsets} of $H$, up to translation by $H$, and $H_A$ denotes the stabilizer in $H$ of such a subset $A$. In particular $H_A=H$ if and only if $A=H$ or $A=\emptyset$, hence if $K$ is a subgroup of $H$, the only term of the form $\gsect{K,K}_H$ in the right hand side of Equation~\ref{sum example} is $\gsect{H,H}_H$, with coefficient 2.\par
On the other hand, by Lemma~\ref{decomposition} again
$$\pi(X\stackrel{\Id}{\to} X)=\sumb{A\subseteq H}{A\,{\rm mod. H}}\gsect{H_A,H_A}_H\mpoint$$
As there are proper non empty subsets $A$ of $H$, there are some terms of the form $\gsect{K,K}_H$ with non zero coefficient in this summation, for some subgroups $K<H$. Hence $\pi(X\stackrel{f}{\to}Y)\neq \pi(X\stackrel{\Id}{\to} X)$, thus $T_U$ does not preserve the defining relations of $\Xi(G)$.
\end{rem}
\npar The following proposition answers Question~(Q1) in the case of left inert bisets~:
\begin{mth}{Proposition} \begin{enumerate}
\item Let $G$ and $H$ be finite groups, and $U$ be a finite left inert $(H,G)$-biset. Then $T_U$ induces a well defined map $t_U:\Xi(G)\to\Xi(H)$, such that $t_U\big(\Gamma(G)\big)\subseteq \Gamma(H)$. Moreover $t_U$ only depends on the isomorphism class of the biset $U$.
\item The map $t_U$ is multiplicative, i.e. $t_U(ab)=t_U(a)t_U(b)$, for any $a,b\in\Xi(G)$. Moreover $t_U(1_{\Xi(G)})=1_{\Xi(H)}$. In particular $t_U$ restricts to group homomorphisms $\Xi(G)^\times\to\Xi(H)^\times$ and $\Gamma(G)^\times\to\Gamma(H)^\times$.
\item If $U$ and $U'$ are finite left inert $(H,G)$-bisets, then $t_{U\sqcup U'}=t_Ut_{U'}$.
\item If $G$, $H$, and $K$ are finite groups, if $U$ is a finite left inert $(H,G)$-biset and $V$ is a finite left inert $(K,H)$-biset, then $t_V\circ t_U=t_{V\times_HU}$. 
\item If $U$ is the identity $(G,G)$-biset, then $t_U$ is the identity map.
\end{enumerate}
\end{mth}
\pf Since $U$ is left inert, Example~\ref{example left inert} shows that $U$ splits as a disjoint union of $(H,G)$-bisets
$$U\cong \bigsqcup_{u\in[U/G]}uG\mpoint$$
It follows that the functor $T_U=\Hom_{G\hbox{-}\mathsf{set}}(U,-)$ is isomorphic to the direct product of the functors $T_{uG}$, for $u\in [U/G]$. Hence to prove Assertion~1, it suffices to consider the case where $U$ is transitive as a right $G$-set.\par
In this case, the functor $T_U$ induces actually a {\em group homomorphism} from $\Xi(G)$ to $\Xi(H)$~: to see this, it suffices to check that the defining relations of $\Xi(G)$ are mapped to relations in $\Xi(H)$. Fix $u\in U$, and denote by $G_u$ its stabilizer in $G$. So $U=uG\cong G/G_u$ as $G$-set, and for any $G$-set $X$, there is a bijection
$$\Hom_{G\hbox{-}\mathsf{set}}(U\op,X)\cong X^{G_u}\mpoint$$
This is actually an isomorphism of $H$-sets if the left action of $H$ on $X^{G_u}$ is defined as follows~: for $h\in H$, there is some $g\in G$ such that $hu=ug$, and the action of $h$ on $X^{G_u}$ is defined by $hx=g^{-1}x$, for $x\in X^{G_u}$.\par
Now let $X\stackrel{f}{\to}Y$ be a morphism of finite $G$-sets such that $X$ splits as a disjoint union $X=X_1\sqcup X_2$ of two $G$-sets. The image $\pi(X\stackrel{f}{\to}Y)$ of this morphism in $\Xi(G)$ is equal to the sum of the images of the morphisms
$$X_1\stackrel{f_1}{\to}f(X_1)\;\;\hbox{and}\;\;X_2\stackrel{f_2}{\to}f(X_2)\mvirg$$
where $f_1$ and $f_2$ are the restrictions of $f$ to $X_1$ and $X_2$, respectively.\par
On the other hand $T_U(X\stackrel{f}{\to}Y)$ is isomorphic to $X_1^{G_u}\sqcup X_2^{G_u}\stackrel{f^{G_u}}{\longrightarrow}Y^{G_u}$ in the category $\mor{H}$, where $f^{G_u}$ is the restriction of $f$ to $X^{G_u}$. The image $\pi\big(T_U(X\stackrel{f}{\to}Y)\big)$ of this morphism in $\Xi(H)$ is equal to the sum of the images of the morphisms
$$X_1^{G_u}\stackrel{f_1^{G_u}}{\longrightarrow}f(X_1^{G_u})\;\;\hbox{and} \;\;X_2^{G_u}\stackrel{f_2^{G_u}}{\longrightarrow}f(X_2^{G_u})\mpoint$$
By Lemma~\ref{image}, since $f(X_1^{G_u})\subseteq f(X_1)^{G_u}$
$$\pi\big(X_1^{G_u}\stackrel{f_1^{G_u}}{\longrightarrow}f(X_1^{G_u})\big)=\pi\big(X_1^{G_u}\stackrel{f_1^{G_u}}{\longrightarrow}f(X_1)^{G_u}\big)\mvirg$$
which is equal to $\pi\left(T_U\big(X_1\stackrel{f_1}{\to}f(X_1)\big)\right)$. It follows that in $\Xi(H)$
$$\pi\big(T_U(X\stackrel{f}{\to}Y)\big)=\pi\left(T_U\big(X_1\stackrel{f_1}{\to}f(X_1)\big)\right)+\pi\left(T_U\big(X_2\stackrel{f_2}{\to}f(X_2)\big)\right)\mvirg$$
hence $T_U$ induces a group homomorphism $t_U:\Xi(G)\to \Xi(H)$. As the functor $T_U$ maps direct product of $G$-sets to direct product of $H$-sets, and the trivial $G$-set to the trivial $H$-set, the morphism $t_U$ is actually a {\em unital ring homomorphism} from $\Xi(G)$ to $\Xi(H)$ (recall that $U$ is assumed transitive as a right $G$-set, here).\par
If $U$ is an arbitrary finite left inert $(H,G)$-biset, then $T_U$ induces the map
$$t_U=\prod_{u\in[U/G]}t_{uG}:\Xi(G)\to \Xi(H)\mpoint$$
It follows from Proposition~\ref{Galois left inert} that $t_U\big(\Gamma(G)\big)\subseteq \Gamma(H)$.\par
Now Assertion~2 follows, since $t_U$ is equal to a product of unital ring homomorphisms. Assertions~3, 4, and 5 are straightforward consequences of the properties of the functor $T_U$.\findemo
\begin{rem}{Remark} It follows that the correspondences $G\mapsto\Xi(G)^\times$ and $G\mapsto\Gamma(G)^\times$ are biset functors for which the biset operations are only defined for {\em left inert bisets}. Equivalently, these correspondences are {\em biset functors without induction}~: the usual basic operations for biset functors are defined for these correspondences, namely {\em restriction} to a subgroup, {\em deflation} from $G$ to a factor group $G/N$ (this is induced by taking {\em fixed points} by $N$ on $G$-sets), transport by isomorphism, and inflation from a factor group. But there is no induction from a subgroup.
\end{rem}
\npar The last result of this Appendix is the computation of the group $\Xi(G)^\times$, when $G$ is abelian~:
\begin{mth}{Theorem} \label{unit group elemab}Let $G$ be a finite abelian group. Then $\Xi(G)^\times$ has an $\F_2$-basis consisting of the elements
$$\left[\begin{array}{l}-\gsect{G,G}_G\\\left.\begin{array}{l}\gsect{G,G}_G-\gsect{S,S}_G\\\gsect{G,G}_G-\gsect{G,S}_G\end{array}\right\}\hbox{for}\;|G:S|=2\mpoint\end{array}\right.$$
In particular
$$\Xi(G)^\times\cong(\F_2)^{2r+1}\mvirg$$
where $r$ is the number of subgroups of index 2 in $G$.
\end{mth}
\pf By Theorem~\ref{units characterization}, the group $\Xi(G)^\times$ is isomorphic to the group of sequences $(u_{T,S})_{(T,S)\in\Pi(G)}$ with values in $\{\pm 1\}$, such that for any $(T,S)\in\Pi(G)$, the map
$$g\in G/S\mapsto u_{\gen{gT},\gen{gS}}/u_{T,S}$$
is a group homomorphism. Switching to an additive notation, the group $\Xi(G)^\times$ is isomorphic to the $\F_2$-vector spaces of sequences $(\lambda_{T,S})_{(T,S)\in \Pi(G)}$, with values in $\F_2$, such that for any $(T,S)\in \Pi(G)$ and any $(g,h)\in G$
$$\lambda_{\gen{ghT},\gen{ghS}}+\lambda_{\gen{gT},\gen{gS}}+\lambda_{\gen{hT},\gen{hS}}+\lambda_{T,S}=0\mpoint$$
If none of $g, h,$ and $gh$ are in $S$, this yields an expression
$$\lambda_{T,S}=\lambda_{\gen{ghT},\gen{ghS}}+\lambda_{\gen{gT},\gen{gS}}+\lambda_{\gen{hT},\gen{hS}}$$
of $\lambda_{T,S}$ as a linear combination of $\lambda_{T',S'}$, for $S'>S$. If $|G:S|>2$, it is always possible to find such elements $g$ and $h$. It follows that the sequence $(\lambda_{T,S})_{(T,S)\in\Pi(G)}$ is entirely determined by the values $\lambda_{T,S}$, for $|G:S|\leq 2$, i.e. the values $\lambda_{G,G}$, $\lambda_{S,S}$ and $\lambda_{G,S}$, where $S$ runs through the set of subgroups of index 2 in $G$. If there are $r$ such subgroups, this gives $2r+1$ such values, so $\dim_{\F_2}\Xi(G)^\times\leq 2r+1$.\par
To prove the converse, and actually the more precise statement in the theorem, it suffices to show that the elements $-1_{\Xi(G)}$, $u_S^G=\gsect{G,G}_G-\gsect{S,S}_G$ and $v_S^G=\gsect{G,G}_G-\gsect{G,S}_G$, for $|G:S|=2$ are linearly independent elements of the $\F_2$-vector space $\Xi(G)^\times$. \par
First as $\gsect{G,G}_G=1_{\Xi(G)}$, and as $\gsect{S,S}_G^2=2\gsect{S,S}$ if $|G:S|=2$, by Proposition~\ref{product formula}, it follows that $(u_S^G)^2=(v_S^G)^2=1$, hence $\{u_S^G,v_S^G\}\subseteq \Xi(G)^\times$. \par
Observe now that for $|G:S|=2$, both elements $u_S^G$ and $v_S^G$ are obtained by inflation from $G/S$ to~$G$ of the corresponding elements $u_\un^{G/S}$ and $v_\un^{G/S}$. The group $C=G/S$ has order 2. The additive group $\Xi(C)$ has a $\Z$-basis consisting of $\gsect{C,C}_C$, $a=\gsect{C,\un}_C$ and $b=\gsect{\un,\un}_C$. The element $\gsect{C,C}_C$ is the identity element $1_{\Xi(C)}$ of $\Xi(C)$, and by Proposition~\ref{product formula}, the products of the other basis elements are given by $a^2=2a$, $b^2=ab=ba=2b$. The group of units $\Xi(C)^\times$ has an $\F_2$-basis consisting of $-1_{\Xi(C)}$, $x=1_{\Xi(C)}-a$ and $y=1_{\Xi(C)}-b$ (actually $\Xi(C)$ is {\em equal} to the group algebra of the multiplicative group $\gen{-1_{\Xi(C)},x,y}\cong (C_2)^3$). Moreover, the set $\{x,y\}$ is an $\F_2$-basis of the kernel $\partial\Xi(G)^\times$ of the deflation map $\Def_{C/C}^C:\Xi(C)^\times\to \Xi(\un)^\times=\{\pm 1_{\Xi(\un)}\}$ (recall that deflation in this case consists in taking fixed points under $C$).\par
Now assume that in $\Xi(G)^\times$, there is a linear relation (with an additive notation) of the form
$$\lambda(-1_{\Xi(G)})+\sum_{|G:S|=2}(\alpha_Su_S^G+\beta_Sv_S^G)=0\mvirg$$
for some coefficients $\lambda, \alpha_S, \beta_S$ in $\F_2$. Fix a subgroup $X$ of index 2 in $G$, and apply $\Def_{G/X}^G$ to this relation. Observe that
$$\Def_{G/X}^Gu_S^G=\Def_{G/X}^G\Inf_{G/S}^Gu_\un^{G/S}=\Inf_{G/SX}^{G/X}\Def_{G/SX}^{G/S}u_\un^{G/S}$$
is equal to 0 if $S\neq X$, and to $u_\un^{G/S}$ if $S=X$. Similarly $\Def_{G/X}^Gv_S^G$ is equal to 0 if $S\neq X$, and to $v_\un^{G/S}$ if $S=X$. It follows that
$$\lambda(-1_{\Xi(G/X)})+\alpha_Xu_\un^{G/X}+\beta_{X}v_\un^{G/X}=0\mvirg$$
hence $\lambda=\alpha_X=\beta_{X}=0$. This completes the proof, since this holds for any~$X$ of index 2 in $G$.\findemo
\begin{rem}{Remark} Using the inclusion $i_G:B(G)\to \Xi(G)$ of Proposition~\ref{from B to Xi}, one can identify $B(G)$ with a subring of $\Xi(G)$. Then it is easy to see that $B(G)^\times$ is the subgroup of $\Xi(G)^\times$ generated by the elements $-\gsect{G,G}_G$ and $\gsect{G,G}_G-\gsect{S,S}_G$, for $|G:S|=2$~: this gives another proof of Matsuda's theorem (\cite{matsuda}), saying that the unit group $B(G)^\times$ of the usual Burnside ring of a finite abelian group $G$ is isomorphic to $(\F_2)^{r+1}$, where $r$ is the number of subgroups of index 2 in $G$.
\end{rem}

Serge Bouc - CNRS-LAMFA, Universit\'e de Picardie, 33 rue St Leu, 80039, Amiens Cedex 01 - France.\\
{\tt email : serge.bouc@u-picardie.fr}\\
{\tt http://www.lamfa.u-picardie.fr/bouc/}
\end{document}
\bibliographystyle{abbrv}
\bibliography{mabib}

\end{document}